\documentclass[
aps,%
pre,%
a4paper,%
final,%
reprint,%
showpacs,%
superscriptaddress,%
nofootinbib
]{revtex4-2}

\usepackage{amsmath}
\usepackage{mathrsfs}  
\usepackage{amssymb}
\usepackage{amsthm}
\usepackage{bm}
\usepackage{multirow}
\usepackage{array}
\usepackage{xcolor}
\colorlet{myrem}{red}

\usepackage[pdftex]{graphicx}
\usepackage[pdftex,%
colorlinks=true,%
bookmarks=true,%
citecolor=blue,%
urlcolor=blue]{hyperref} 
\usepackage[caption=false]{subfig}

\newcommand{\field}[1]{\mathbb{#1}}
\newcommand{\OP}[1]{\mathscr{#1}}

\newcommand{\tp}{\intercal}
\DeclareMathOperator*{\supp}{supp}
\let\Re\relax
\DeclareMathOperator{\Re}{Re}
\let\Im\relax
\DeclareMathOperator{\Im}{Im}

\DeclareMathOperator{\Span}{Span}
\newcommand{\melem}[1]{\mathfrak{#1}}
\newcommand{\what}[1]{\widehat{#1}}

\graphicspath{{./figs/}}
\DeclareGraphicsExtensions{.pdf,.jpeg,.png}
\newcommand{\wtilde}[1]{\widetilde{#1}}

\newcommand{\fs}[1]{\mathsf{#1}}

\DeclareMathOperator{\diag}{diag}

\newcommand{\ovl}[1]{\overline{#1}}

\let\Re\relax
\DeclareMathOperator{\Re}{Re}
\let\Im\relax
\DeclareMathOperator{\Im}{Im}
\newcommand{\vv}[1]{\boldsymbol{#1}}
\newcommand{\vs}[1]{\boldsymbol{#1}}

\DeclareMathOperator{\sgn}{sgn}

\newtheorem{rem}{Remark}
\usepackage{enumitem}

\begin{document}
\title{Nonreflecting Boundary Condition for the free 
Schr\"{o}dinger equation for hyperrectangular computational domains}

\author{Samardhi Yadav}
\affiliation{Department of Physics, Indian Institute of Technology Delhi, Hauz
Khas, New Delhi–110016, India}%
\email{Samardhi@physics.iitd.ac.in}
\author{Vishal Vaibhav}
\affiliation{Department of Physics, Indian Institute of Technology Delhi, Hauz
Khas, New Delhi–110016, India}%
\affiliation{Optics and Photonics Center, Indian Institute of Technology Delhi,
Hauz Khas, New Delhi–110016, India}%
\email{vishal.vaibhav@gmail.com}
\date{\today}

\begin{abstract}
In this article, we discuss efficient ways of implementing the transparent boundary 
condition (TBC) and its various approximations for the free Schr\"{o}dinger equation 
on a hyperrectangular computational domain in $\field{R}^d$ with periodic boundary conditions
along the $(d-1)$ unbounded directions. In particular, we consider Pad\'e approximant 
based rational approximation of the exact TBC and a spatially local form of the exact 
TBC obtained under its high-frequency approximation. 
For the spatial discretization, we use a Legendre-Galerkin spectral method 
with a boundary-adapted basis to ensure the bandedness of the resulting linear system. 
Temporal discretization is then addressed with two one-step methods, namely, the 
backward-differentiation formula of order 1 (BDF1) and the trapezoidal rule (TR). 
Finally, several numerical tests are presented to demonstrate the effectiveness of the 
methods where we study the stability and convergence behaviour empirically.
\end{abstract}


\maketitle
\section*{Notations}
The set of non-zero positive real numbers ($\field{R}$) is denoted by
$\field{R}_+$. For any complex number $\zeta$, $\Re(\zeta)$ and $\Im(\zeta)$ refer to the real
and the imaginary parts of $\zeta$, respectively.
\section{Introduction}
This article is an extension of the work presented in~\cite{SV2023} where we considered 
the numerical implementation of transparent boundary conditions for the free Schr\"{o}dinger 
equation on a rectangular computational domain with periodic boundary condition along one 
of the unbounded directions. In this work, we consider the elliptic as well as
the non-elliptic Schr\"odinger equation in $\field{R}^d$, $d=2,3$, with periodic boundary
condition in $d-1$ dimensions. Both of these models can be encountered in
optics~\cite{MN1990}. For instance, in the case of an uniform isotropic media,
the propagation of optical pulses under slow varying envelope approximation works out to be
\begin{equation}\label{eq:model-ivp}
\left[i\left(\partial_{x_3}+\frac{1}{v_g}\partial_{\tau}\right)-\frac{1}{2k}\nabla^2_{\perp}
+\frac{k''}{2}\partial^2_{\tau}\right]A(\vv{x},\tau)=0
\end{equation}
where $A(\vv{x},\tau)$ is the envelop, $k=k(\omega_0)$ is the wavenumber corresponding to 
$\omega_0$, the central frequency. The quantity $v_g=[k'(\omega_0)]^{-1}$ is the
group velocity and $k''(\omega_0)$ is related to group velocity dispersion. Depending on the 
sign of the coefficients, the equation above can turn out to be elliptic or non-elliptic. Let us 
emphasize that this work can also be considered as a preparation for treating more general 
models which include variable potentials and nonlinearity within the 
pseudo-differential approach on the lines of~\cite{V2019}.

For the purpose of this exposition, we consider the initial value problem (IVP)
corresponding to the Schr\"odinger equation in the following standard form:
\begin{equation}\label{eq:2D-SE}
i\partial_tu+\partial^2_{x_1}u+\beta\nabla_{\perp}u=0,
\quad(\vv{x},t)\in\field{R}^d\times\field{R}_+,
\end{equation}
where $\beta=\pm1$ and the initial condition $u(\vv{x},0)=u_0(\vv{x})$ is assumed to be compactly 
supported in the computational domain, $\Omega_i$. We are interested
in the efficient numerical solution of~\eqref{eq:2D-SE} on a hyperrectangular 
computational domain in $\field{R}^d$ with periodic boundary conditions along 
$x_j,\,j>1$.

The primary objective of this work is generalization of the novel Pad\'e approach 
previously applied to a rectangular domain in our work~\cite{YV2024I}, to a hyperrectangular domain, 
leveraging its advantages of time-step-independent computational complexity. 

The transparent boundary operator for the Schr\"odinger equation on a 
hyperrectangular domain has the form $\sqrt{\partial_t-i\triangle_{\Gamma}}$ where 
$\triangle_{\Gamma}$ is the Laplace-Beltrami operator. Note that these operators are 
nonlocal, both in time as well as space. The earliest numerical implementation 
for the operator $\sqrt{\partial_t-i\triangle_{\Gamma}}$ was given by 
Menza~\cite{Menza1997} who introduced a Pad\'e approximant based 
rational approximation for the complex square root function to obtain effectively 
local boundary conditions. However, this method failed to address the problems
arising at the corners of the rectangular domain. 
In our recent work~\cite{YV2024I}, we reported a novel Pad\'e approach that
effectively localizes the transparent boundary operator, offering a 
computationally efficient solution that handles corners adequately. 
In this paper, we developed an extension of this approach for the 
hyperrectangular domain, which reduces the nonlocality in space to solving
two-dimensional Schr\"odinger equations satisfied by certain auxiliary
functions on the boundary faces of the hyperectangular computational domain,
while the temporal nonlocality is handled through rational approximation for 
$\sqrt{\partial_t}$ operator acting on the auxiliary
functions. The novel-Pad\'e method developed in our work ensures that the 
computational cost and the memory requirement does not grow with increasing
time-steps. Furthermore, we provide a detailed comparison of various
approximations for the TBCs, including convolution quadrature (CQ) approach, conventional Pad\'e approach 
(CP), novel Pad\'e approach (NP), and high-frequency approximation. 
Note that CQ scheme remains nonlocal while others are effectively local.

For the space discretization, we use a Legendre-Galerkin method where a boundary-adapted 
basis is constructed to ensure the bandedness of the resulting linear system.
At the discrete level, the TBCs are formulated as a Robin-type boundary condition 
at each time-step. In order to homogenize these Robin-type maps needed for the 
construction of the basis, we also discuss a boundary-lifting process. 
Temporal discretization is then addressed with two one-step time marching methods, 
namely, BDF1 and TR. Let us emphasize that an unsuitable temporal discretization of the
nonlocal operators may destroy the overall stability of the numerical 
scheme~\cite{Mayfield1989,PM2015}. To avoid this, we try to match the time-stepping 
methods in the discretization of the boundary conditions to that of the interior
problem.

The paper is organized in the following manner: Sec.~\ref{sec:tbcs-2d} presents the
derivation of TBCs followed by a discussion of numerically tractable forms of
the TBCs for the 2D problem. This section also discusses the high-frequency 
approximation of the TBCs. Sec.~\ref{sec:tbcs-3d} presents the derivation of
TBCs followed by discussion of novel Pad\'e based representation of TBCs for the
3D problem. 
In Sec.~\ref{sec:numerical-implementation-2d} and
Sec.~\ref{sec:numerical-implementation-3d}, we discuss the time-discretization of 
TBCs and the numerical solution of the IBVP using Legendre-Galerkin method 
with a boundary-adapted basis for the 2D and 3D problems, respectively. We then
confirm the efficiency and stability of our numerical schemes with several numerical 
tests presented in Sec.~\ref{sec:numerical-experiments-2d} and
Sec.~\ref{sec:numerical-experiments-3d}. Finally, we conclude this paper in 
Sec.~\ref{sec:conclusion}.


\section{TBC--2D}\label{sec:tbcs-2d}
\begin{figure}[!htbp]
\begin{center}
\def\myscale{1}
\includegraphics[scale=\myscale]{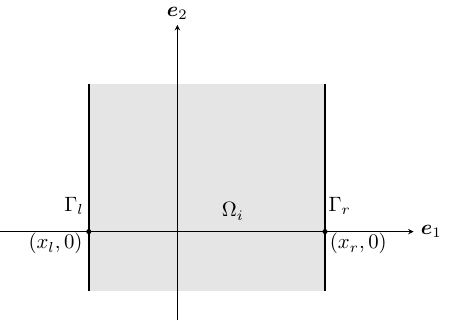}
\end{center}
\caption{\label{fig:rect-domain} The figure shows a rectangular domain with
periodic boundary conditions along vertical axis.}
\end{figure}

Consider a rectangular computational domain with boundary segments parallel to
one of the axes (see Fig.~\ref{fig:rect-domain}) with $\Omega_i =(x_l,x_r)\times[-d,d)$
also termed as \emph{interior} domain. By using the decomposition of 
$u(\vv{x},t)\in \fs{L}^2(\field{R}\times[-d,d))=\fs{L}^2(\Omega_i)
\oplus \fs{L}^2(\Omega_e)$ where $\Omega_{e}=(\field{R}\times[-d,d))\setminus\overline{\Omega}_i$ 
also termed as \emph{exterior} domain. With this in mind, the equivalent formulation of 
problem~\eqref{eq:2D-SE} in $\field{R}^2$ can be stated as  
\begin{equation}\label{eq:split-IVP}
\begin{split}
\text{interior}:\;
&\left\{\begin{aligned}
&i\partial_tv+ \triangle v=0,\quad(\vv{x},t)\in\Omega_i\times\field{R}_+,\\
&v(\vv{x},0)= u_0(\vv{x}),\quad\vv{x}\in\Omega_i,\\
\end{aligned}\right.\\
\text{exterior}:\;
&\left\{\begin{aligned}
&i\partial_tw+\triangle w=0,\quad(\vv{x},t)\in\Omega_e\times\field{R}_+,\\
&w(\vv{x},0)=0,\quad\vv{x}\in\Omega_e,\\
&\lim_{|{x_1}|\rightarrow\infty} |w(\vv{x},t)|=0.
\end{aligned}\right.\\
\text{continuity}:\;
&v|_{\Gamma}=w|_{\Gamma},\;
\partial_{n}v|_{\Gamma}=\partial_{n}w|_{\Gamma},
\end{split}
\end{equation}
where $\triangle=(\partial_{x_1}^2+\beta\partial_{x_2}^2)$ and $\Gamma$ denotes the common boundary of 
the interior and the exterior domains. Also, the initial data is assumed to be compactly supported in $\Omega_i$.
For the sake of completeness, we reproduce the results from~\cite{SV2023,V2019}
on the derivation of the TBCs for the system under consideration. Assuming the 
periodicity behaviour along $x_2$, we can express the field $w(x_1,x_2,t)$ in terms 
of a Fourier series as:
\begin{equation}
\begin{split}
&w(x_1,x_2,t)=\sum_{m\in\field{Z}}\widetilde{w}_m(x_1,t)e^{i\zeta_m x_2},\\
&W(x_1,x_2,z)=\sum_{m\in\field{Z}}\widetilde{W}_m(x_1,z)e^{i\zeta_m x_2},
\end{split}
\end{equation}
where $\widetilde{W}_m(x_1,z)=\OP{L}_t[\widetilde{w}_m(x_1,t)]$ denote the Laplace
transform of the field.
By employing the Fourier expansion with respect to $x_2$
and then taking the Laplace transform of the right exterior problem yields
\begin{equation}
(\partial_{x_1}^2+\alpha_m^2)\wtilde{W}_m(x_1,z)=0,\quad x_1\in(x_r,\infty),
\end{equation}
where $\alpha_m=\sqrt{(iz-\beta\zeta_m^2)}$ such that $\sqrt{\cdot}$ denotes the branch
with $\Im(\alpha_m)>0$. The bounded solution works out to be 
\begin{equation}
\partial_{x_1}\wtilde{W}_m(x_1,z) =i\alpha_m\wtilde{W}(x_1,z).
\end{equation}
To facilitate the inverse Laplace transform, we rewrite the above equation as
\begin{equation}\label{eq:dtn-map-freq}
\OP{L}^{-1}[\partial_{x_1}\wtilde{W}_m]
=\OP{L}^{-1}[i\alpha_m^{-1}]\star\OP{L}^{-1}[\alpha_m^2\wtilde{W}_m],
\end{equation}
where `$\star$' represents the convolution operation.
In order to carry out the inverse Laplace transformation, let us observe that 
$\OP{L}^{-1}\left[{1}/{\sqrt{z}}\right](t)=H(t)/{\sqrt{\pi t}}$ where $H(t)$
denotes the Heaviside step-function. In the following, we use 
the fractional operators which are denoted by 
$\partial_t^{\alpha},\;\alpha\in\field{R}$ (see Miller and Ross~\cite{MR1993}).
For $\alpha<0$, we obtain the Riemann-Liouville fractional integral of order 
$|\alpha|$ while $\alpha>0$ corresponds to fractional derivatives which are 
introduced as integral order derivatives of fractional integrals. Using the shifting 
property of Laplace transforms and introducing fractional operators of order $1/2$, we 
obtain the following result
\begin{equation}\label{eq:dtn-map}
\begin{split}
\partial_{x_1}\wtilde{w}_m(x_1,t)
&=e^{i\pi/4}e^{-i \beta\zeta_m^2t}\partial_t^{-1/2}e^{i\beta\zeta_m^2t}\times\\
&\quad\left[i\partial_t\wtilde{w}_m(x_1,t)-\beta\zeta_m^2\wtilde{w}_m(x_1,t)\right].
\end{split}
\end{equation}
In order to obtain the inverse Fourier transform of the above relation, let us 
introduce the kernel which is defined in a distributional sense as
\begin{equation}\label{eq:kernel-periodic}
\mathcal{H}(x_2,t)=\frac{1}{2d}\sum_{m\in\field{Z}}e^{i(\zeta_m x_2-\beta\zeta_m^2t)}.
\end{equation}
Next by introducing the integral definition for the operator $\partial_t^{-1/2}$ and 
taking the inverse Fourier transform, the DtN map reads as
\begin{widetext}
\begin{equation}\label{eq:IPS-TBC}
\begin{split}
\partial_{x_1}{w}(\vv{x},t)
&=\frac{e^{i\pi/4}}{\sqrt{\pi}}\int_0^t\int_{\field{R}}
\left[i\partial_{\tau}{w}(x_1,x'_2,\tau)+\beta\partial^2_{x'_2}w(x_1,x'_2,\tau)
\right]
\frac{\mathcal{H}(x_2-x'_2,t-\tau)}{\sqrt{t-\tau}}dx_2'd\tau\\
&=-(\partial_{t}-i\beta\partial^2_{x_2})\frac{e^{-i\pi/4}}{\sqrt{\pi}}
\int_0^t\int_{\field{R}}w(x_1,x'_2,\tau)
\frac{\mathcal{H}(x_2-x'_2,t-\tau)}{\sqrt{t-\tau}}dx_2'd\tau,
\end{split}
\end{equation}
\end{widetext}
where the last step is obtained using the following result from
the calculus of fractional operators~\cite{MR1993}:
\begin{multline}\label{eq:ISP-Op2}
\partial_t\left[e^{-i\beta\zeta_m^2t}\partial_t^{-1/2}e^{i\beta\zeta_m^2t}
\wtilde{w}_m(x_1,t)\right]\\
=e^{-i\beta\zeta_m^2t}\partial_t^{-1/2}
e^{i\beta\zeta_m^2t}\partial_t \wtilde{w}_m(x_1,t).
\end{multline}
In order to express the DtN map compactly, we introduce the notation 
\begin{multline}\label{eq:Op-ISP-TBC}
(\partial_t-i\beta\partial_{x_2}^2)^{-1/2}f({x}_2,t)\\
=\frac{1}{\sqrt{\pi}}\int_0^t\int_{\field{R}}f(x_1,x'_2,\tau)
\frac{\mathcal{H}(x_2-x'_2,t-\tau)}{\sqrt{t-\tau}}dx_2'd\tau,
\end{multline}
so that
\[ (\partial_t-i\beta\partial_{x_2}^2)^{1/2}f
=(\partial_t-i\beta\partial_{x_2}^2)[(\partial_t-i\beta\partial_{x_2}^2)^{-1/2}f].\]
By using the continuity relations defined in~\eqref{eq:split-IVP}, the DtN
map on the right segment can be compactly written as
\begin{equation}
\partial_{x_1}{v}(\vv{x},t)
+e^{-i\pi/4}(\partial_t-i\beta\partial^2_{x_2})^{1/2}v(\vv{x},t)=0.
\end{equation}
An equivalent formulation of the IVP~\eqref{eq:2D-SE} on the computational
domain $\Omega_i = (x_l,x_r)\times[-d,d)$ with periodic boundary condition
along the axis $\vv{e}_2$ is given by
\begin{equation}\label{eq:2D-SE-CT}
\left\{\begin{aligned}
&i\partial_tu+(\partial_{x_1}^2+\beta\partial_{x_2}^2)u=0,\quad (\vv{x},t)\in\Omega_i\times\field{R}_+,\\
&u(\vv{x},0)=u_0(\vv{x})\in \fs{L}^2(\Omega_i),
\quad\supp\,u_0\subset\Omega_i;\\
&\partial_{n_1}{u}+
e^{-i\pi/4}(\partial_{t}-i\beta\partial^2_{x_2})^{1/2}u=0,
\quad\vv{x}\in\Gamma_l\cup\Gamma_r,\\
&u(x_1, x_2+2d,t)=u(x_1, x_2,t),\quad t>0.
\end{aligned}\right.
\end{equation}
Note that the operators present in the 
DtN maps are non-local in space as well 
time. For an efficient numerical implementation of this IBVP, we first focus
on writing a numerically tractable form of these non-local operators at 
a continuous level in Sec.~\ref{sec:CT-CQ}. We revisit the novel Pad\'e based 
approach from~\cite{YV2024I} in Sec.~\ref{sec:CT-NP} which is then contrasted 
with existing Pad\'e approach for such operators in Sec.~\ref{sec:CT-CP}. 
Finally, we present an approximate form of these operators under high-frequency 
approximation from~\cite{YV2024II} in Sec.~\ref{sec:CT-HF}. 
\subsection{Fractional operators based approach}\label{sec:CT-CQ}
\begin{figure}[!htbp]
\begin{center}
\def\myscale{1}
\includegraphics[scale=0.7]{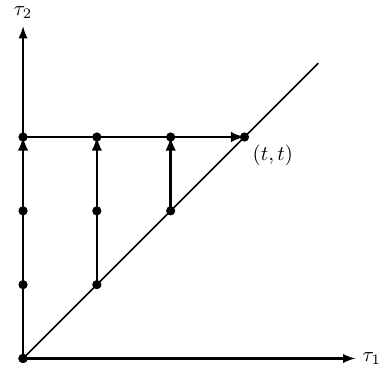}
\end{center}
\caption{\label{fig:IVP-auxi-periodic} A schematic is shown in this figure which
depicts the evolution of the auxiliary field $\varphi(x_1,x_2,\tau_1,\tau_2)$ 
in the $(\tau_1,\tau_2)$-plane starting from the diagonal which also serves as 
initial conditions for solving the IVPs (arrow in the line depicts the 
evolution direction in time). Note that TBCs on $\partial{\Omega_i}$ require 
the history of the auxiliary field along the horizontal line up to the diagonal in the 
$(\tau_1,\tau_2)$-plane.}
\end{figure}

In this section, we explore the possibility of writing the DtN maps in terms 
of (time-)fractional operators. To achieve this, let us recall the form of the DtN map given 
by~\eqref{eq:dtn-map} in terms of the Fourier variables so that
\begin{equation}
\begin{split}
&\partial_{x_1}\wtilde{w}_m(x_1,t)\\
&=e^{i\frac{\pi}{4}}e^{-i \beta\zeta_m^2t}\partial_t^{-1/2}e^{i\beta\zeta_m^2t}
\left[i\partial_t -\beta\zeta_m^2\right]\wtilde{w}_m(x_1,t)\\
&= ie^{i\frac{\pi}{4}}e^{-i\beta\zeta_m^2t}\partial_t^{1/2}e^{i\beta\zeta_m^2t} \wtilde{w}_m(x_1,t).
\end{split}
\end{equation}
Taking inverse Fourier transform yields
\begin{equation*}
\begin{split}
\partial_{n_1}{w}
&=-e^{i\frac{\pi}{4}}\sum_{m\in\field{Z}}e^{-i\beta\zeta_m^2t}\partial^{1/2}_t
\left[e^{i\beta\zeta_m^2t}\widetilde{w}_m(x_1,t)e^{i\zeta_mx_2}\right]\\
&=\left.-e^{i\frac{\pi}{4}}\partial^{1/2}_{t'}\sum_{m\in\field{Z}}
\left[e^{-i\beta\zeta_m^2(t-t')}\widetilde{w}_m(x_1,t')e^{i\zeta_m x_2}\right]\right|_{t'=t}.
\end{split}
\end{equation*}
Introducing the auxiliary function $\varphi(x_1,x_2,\tau_1,\tau_2)$ such that
\begin{equation}
\begin{split}
&\varphi(x_1,x_2,\tau_1,\tau_2)\\
&=\sum_{m\in\field{Z}}e^{-i\beta\zeta_m^2(\tau_2-\tau_1)}\widetilde{w}_m(x_1,\tau_1)
e^{i\zeta_mx_2}\\
&=\sum_{m\in\field{Z}}\left[e^{i\beta\zeta_m^2\tau_1}\widetilde{w}_m(x_1,\tau_1)\right]
e^{-i\beta\zeta_m^2\tau_2+i\zeta_mx_2}.
\end{split}
\end{equation}
The DtN map on segment $\Gamma_r$ can now be expressed compactly as
\begin{equation}
\begin{split}
&\partial_{x_1}{w}(\vv{x},t)=\left.-e^{i\frac{\pi}{4}}
\partial^{1/2}_{\tau_1}\varphi(x_1,x_2,\tau_1,\tau_2)\right|_{\tau_1,\tau_2=t}.
\end{split}
\end{equation}
By definition, the numerical implementation of this fractional derivative requires 
the history of the auxiliary function from the start of the computation which can 
be facilitated by the IVP satisfied by the auxiliary function: 
\begin{equation}\label{eq:FSE-auxi-ivp}
\left\{\begin{aligned}
&(\tau_2,x_2)\in(\tau_1,t]\times\field{R}:\\
&\quad [i\partial_{\tau_2}+\beta\partial^2_{x_2}]\varphi(x_1,x_2,\tau_1,\tau_2)=0,\\
&\quad \varphi(x_1,x_2+2d,\tau_1,\tau_2)=\varphi(x_1,x_2,\tau_1,\tau_2),\\
&\quad \varphi(x_1,x_2,\tau_1,\tau_1)=w(x_1,x_2,\tau_1).
\end{aligned}\right.
\end{equation}
Note that the history is needed for all $\tau_1\in[0,t]$ which means that we must 
solve the IVP for all $\tau_1\in[0,t]$ over $\Gamma_r$. It is interesting to note 
that the IVP satisfied by the auxiliary function on the segment $\Gamma_r$ turns 
out to be a one-dimensional Schr\"{o}dinger equation. 
The computation of required history can be understood with the help of the
Fig.~\ref{fig:IVP-auxi-periodic} where we note that the 
history of the field is needed along the horizontal line up to the diagonal in the
 $(\tau_1,\tau_2)$-plane. The arrows point towards the direction of evolution
of the auxiliary function starting from the diagonal values which also serve as the 
initial conditions. 

The DtN maps present in~\eqref{eq:2D-SE-CT} in terms of a time-fractional 
derivative operator reads as
\begin{equation}\label{eq:dtn-maps-cq}
\partial_{n_1}{u} = - e^{-i\frac{\pi}{4}}\left.\partial^{1/2}_{\tau_1}
\varphi(x_1,x_2,\tau_1,\tau_2)\right|_{\tau_1=\tau_2=t}.
\end{equation}
Here, we would like to mention that fractional operator based approach for the
numerical implementation of DtN maps is amenable to convolution quadratures but
this approach requires the history from the start of the computation that grows 
with each time-step making it expensive computationally as well as from a 
storage point of view.
\subsection{Novel Pad\'e based approach}\label{sec:CT-NP}
In this section, we would like to explore the Pad\'e approximant based representation 
for the $1/2$-order temporal derivative in~\eqref{eq:dtn-maps-cq} and its
other ramifications. This would allow us to obtain an effectively local approximation 
for the fractional derivative operator. We contrast this with the situation where 
the history grows with each time-step making it expensive computationally as well 
as from a storage point of view. 

\begin{figure}[!htbp]
\begin{center}
\includegraphics[scale=0.7]{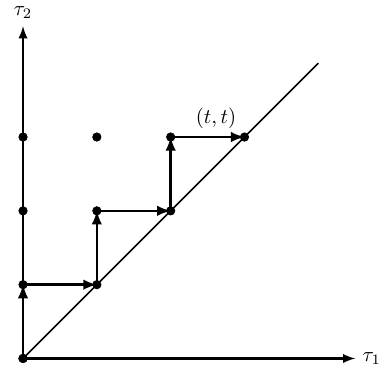}
\end{center}
\caption{\label{fig:IVP-auxi-pade}A schematic depiction of the evolution of the 
auxiliary fields $\varphi_{k,a_1}(x_2,\tau_1,\tau_2)$ in the $(\tau_1,\tau_2)$-plane 
is provided in this figure. }
\end{figure}

The diagonal Pad\'{e} approximants for $\sqrt{1-z}$, with negative real axis as
the branch-cut, can be used to obtain rational approximation $R_M(z)$ for the 
function $\sqrt{z}$ as
\begin{equation}\label{eq:sqrt-pade}
\begin{split}
&R_M(z)
=b_0-\sum_{k=1}^M\frac{b_k}{z+\eta^2_k},\\
&\text{where}\quad
\left\{\begin{aligned}
& b_0=2M+1,\\
&b_k = \frac{2\eta^2_k(1+\eta^2_k)}{2M+1},\;\eta_k = \tan\theta_k,\\
&\theta_k =\frac{k\pi}{2M+1},\;k=1,2,\ldots,M.
\end{aligned}\right.
\end{split}
\end{equation}
Here, the expression for the diagonal Pad\'e approximant is taken from~\cite{V2011}. 
Let us start by introducing a compact notation for the auxiliary functions as:
\begin{equation}
\varphi_{a_1}(x_2,\tau_1,\tau_2)=\varphi(x_{a_1},x_2,\tau_1,\tau_2),\quad a_1\in\{l,r\}.
\end{equation}
Let $z\in\field{C}$, then the Laplace transform of the auxiliary function can be denoted as  
\begin{equation}
\OP{L}_{\tau_1}[\varphi_{a_1}(x_2,\tau_1,\tau_2)]=\varPhi_{a_1}(x_2,z,\tau_2).
\end{equation}
The DtN maps described in~\eqref{eq:dtn-maps-cq} on the segments $\Gamma_{a_1}$  
can be restated as 
\begin{equation}
\partial_{x_1}{u}(\vv{x},t)\pm e^{-i\frac{\pi}{4}}\left.\OP{L}^{-1}_{z}\left[
\sqrt{z}\varPhi_{a_1}\right]\right|_{\tau_1=\tau_2=t}=0.
\end{equation}
Inserting the rational approximation of order $M$ for the function $\sqrt{z}$ as
\begin{equation}\label{eq:maps-npade1}
\partial_{x_1}{u}(\vv{x},t)\pm e^{-i\frac{\pi}{4}}\left.\OP{L}^{-1}_{z}\left[
R_M(z)\varPhi_{a_1}\right]\right|_{\tau_1=\tau_2=t}\approx 0.
\end{equation}
The map expressed above is still non-local, however, the rational nature of the 
expression allows us to obtain an effectively local approximation for the DtN 
operation by introducing the auxiliary fields $\varphi_{k,a_1}(x_2,\tau_1,\tau_2)$ 
at the boundary segment $\Gamma_{a_1}$ as
\begin{equation}
(z+\eta^2_k)^{-1}\varPhi_{a_1}(x_1,x_2,z,\tau_2)=\varPhi_{k,a_1}(x_2,z,\tau_2),
\end{equation}
where $ k=1,2,\ldots,M$. In the physical space, every $\varphi_{k,a_1}$ satisfy the IVP given by
\begin{equation}\label{eq:auxi-npade}
(\partial_{\tau_1}+\eta^2_k)\varphi_{k,a_1}(x_2,\tau_1,\tau_2)
=\varphi(x_1,x_2,\tau_1,\tau_2),
\end{equation}
with the initial conditions set to $\varphi_{k,a_1}(x_2,0,\tau_2)=0$. The DtN map 
now reads as 
\begin{equation}\label{eq:maps-npade}
\partial_{x_1}{u}(\vv{x},t)+e^{-i\frac{\pi}{4}}\left[ b_0 u(\vv{x},t) -
\sum_{k=1}^M b_k\varphi_{k,a_1}(x_2,t,t)\right]=0.
\end{equation}
The solution to the ODEs~\eqref{eq:auxi-npade} work out to be
\begin{equation}
\varphi_{k,a_1}(x_2,\tau_1,\tau_2)
=\int_0^{\tau_1}e^{-\eta_k^2(\tau_1-s)}\varphi(x_{a_1},x_2,s,\tau_2)ds.
\end{equation}
It is easy to verify that the auxiliary fields $\phi_{k,a_1}(x_2,\tau_1,\tau_2)$ 
satisfy the IVPs given by
\begin{equation}\label{eq:ivp-pade-auxi}
[i\partial_{\tau_2}+\beta\partial^2_{x_2}]\varphi_{k,a_1}(x_2,\tau_1,\tau_2)=0.
\end{equation}
The evolution of the auxiliary fields $\varphi_{k,a_1}(x_2,\tau_1,\tau_2)$
in the $(\tau_1,\tau_2)$-plane can be understood with the help of the schematic 
shown in Fig.~\ref{fig:IVP-auxi-pade}. The diagonal to diagonal computation 
of the fields $\varphi_{k,a_1}(x_2,\tau_1,\tau_2)$ consists of first advancing 
the fields using the IVP established in~\eqref{eq:ivp-pade-auxi} and then using 
the ODEs in~\eqref{eq:auxi-npade} for the second movement. The novel Pad\'e 
approach makes the numerical scheme efficient from a storage point of view 
which is also clear from the schematic presented. 

\subsection{Conventional Pad\'{e} based approach}\label{sec:CT-CP}
In this section, we focus on obtaining an effectively local form for the full operator
$(\partial_{t}-i\beta\partial^2_{x_2})^{1/2}$ present in the DtN map in~\eqref{eq:2D-SE-CT}
using Pad\'e approximants. One can obtain the rational approximation for the 
function $\sqrt{Z}$ where $ Z=z+i\beta\zeta_m^2$ using~\eqref{eq:sqrt-pade} now in terms of $Z$.
The DtN map described in~\eqref{eq:2D-SE-CT} on the segments $\Gamma_{a_1}$  
can be restated as 
\begin{equation}
\partial_{x_1}{u}(\vv{x},t)\pm e^{-i\frac{\pi}{4}}\OP{F}^{-1}_{\zeta_m}\OP{L}^{-1}_{z}\left[
\sqrt{Z}U_m(x_{a_1},z)\right]=0.
\end{equation}
Inserting the rational approximation of order $M$ for the function $\sqrt{Z}$ as
\begin{equation}\label{eq:maps-cpade1}
\partial_{x_1}{u}(\vv{x},t)\pm e^{-i\frac{\pi}{4}}\OP{F}^{-1}_{\zeta_m}\OP{L}^{-1}_{z}\left[
R_M(Z)U_m\right]\approx 0.
\end{equation}
Let $k=1,2,\ldots,M$, then introducing the auxiliary fields 
$\{\varphi_k(x_2,t)\leftrightarrow\wtilde{\varPhi}_{k,m}(z),\;m\in\field{Z},\,z\in\field{C}\}$ 
at the boundary segments $\Gamma_{a_1}$ as
\begin{equation}
(z+i\beta\zeta_m^2+\eta^2_k)^{-1}\wtilde{U}_m(x_1,z)=\wtilde{\varPhi}_{k,m}(z).
\end{equation}
In the physical space, every $\varphi_k$ satisfy the IVP given by
\begin{equation}\label{eq:auxi-cpade}
(i\partial_t +\beta \partial^2_{x_2} +i\eta^2_k)\varphi_k(x_2,t)= iu(\vv{x},t),
\end{equation}
with the initial conditions set to $\varphi_k(x_2,0)=0$. The DtN 
map~\eqref{eq:maps-cpade1} now reads as 
\begin{equation}\label{eq:maps-cpade}
\partial_{x_1}u(\vv{x},t)\pm e^{-i\pi/4}\left[ b_0 u(\vv{x},t) -
\sum_{k=1}^M b_k\varphi_k(x_2,t)\right]=0.
\end{equation}
\subsection{High-Frequency Approximation}\label{sec:CT-HF}
In this section, we focus on obtaining the absorbing boundary conditions for the 
operator $(\partial_{t}-i\beta\partial^2_{x_2})^{1/2}$ under high-frequency approximation.
The high-frequency approximation allows us to make the TBCs local with respect to
spatial variables. To this end, we consider $m$-th Fourier component as  
\begin{equation}
\OP{L}^{-1}[i\alpha_m^{-1}]=\frac{1}{2\pi}\int(iz-\beta\zeta^2_m)^{-{1}/{2}}e^{zt}dz.
\end{equation}
Setting $\xi=zt$, we have
\begin{equation*}
\begin{split}
&\OP{L}^{-1}\left[\frac{i}{\alpha_m}\right]
=\frac{t^{-\frac{1}{2}}}{2\pi}\int_{a+i\field{R}}\frac{e^{\xi}}{\sqrt{i\xi-\beta\zeta^2_mt}}d\xi\\
&=\frac{t^{-\frac{1}{2}}}{2\pi}\int_{a+i\field{R}}
\left(1-\frac{\beta\zeta^2_mt}{i\xi}\right)^{-\frac{1}{2}}\frac{e^{\xi}d\xi}{\sqrt{i\xi}}\quad(a>0)\\
&\sim\frac{t^{-\frac{1}{2}}}{2\pi}\int_{a+i\field{R}}
\left(\frac{1}{(i\xi)^{\frac{1}{2}}}
+\frac{\beta\zeta^2_mt}{2(i\xi)^{\frac{3}{2}}}
+\frac{3\zeta^4_mt^2}{8(i\xi)^{\frac{5}{2}}}
+\ldots\right)e^{\xi}d\xi\\
&\sim \frac{e^{i\pi/4}}{\Gamma(\tfrac{1}{2})}t^{-\frac{1}{2}}
+\frac{e^{-i\pi/4}}{2\Gamma(\tfrac{3}{2})}t^{\frac{1}{2}}\beta\zeta_m^2
-\frac{3e^{i\pi/4}}{8\Gamma(\tfrac{5}{2})}t^{\frac{3}{2}}\beta^2\zeta_m^4
+\ldots
\end{split}
\end{equation*}
Taking the inverse Fourier transform of~\eqref{eq:dtn-map-freq}, we obtain the following form of
DtN map as a result of high-frequency approximation on $\Gamma_r$:
\begin{equation}\label{eq:dtn-maps-hf}
\partial_{x_1}u+e^{-i\frac{\pi}{4}}\partial_t^{\frac{1}{2}}u
-\frac{1}{2}e^{i\frac{\pi}{4}}\beta\partial_{x_2}^2\partial_t^{-\frac{1}{2}}
u=0\mod{(\partial_t^{-\frac{3}{2}})}.
\end{equation}
\section{TBC--3D}\label{sec:tbcs-3d}
In this section, we present an extension of the ideas from $\field{R}^2$ to 
$\field{R}^3$ allowing periodic boundary conditions to be imposed along two of 
the axes. As it turns out, extension of all the operators defined in 
Sec.~\ref{sec:tbcs-2d} can be carried out in a natural way.

\begin{figure}[!tbp]
\begin{center}
\def\myscale{1}
\includegraphics[scale=\myscale]{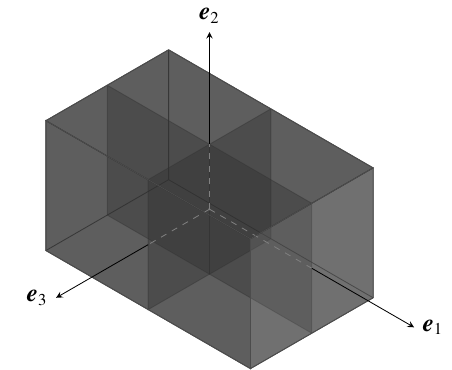}
\end{center}
\caption{\label{fig:cuboidal} The figure shows a rectangular cuboidal domain with 
boundary faces parallel to one of the axes.}
\end{figure}

The initial value problem (IVP) corresponding to the free Schr\"{o}dinger equation 
formulated on $\field{R}^3$ reads as 
\begin{equation}\label{eq:3D-SE}
i\partial_tu+ (\partial_{x_1}^2+\beta\partial_{x_2}^2+\beta\partial_{x_3}^2)u=0,
\quad(\vv{x},t)\in\field{R}^3\times\field{R}_+,
\end{equation}
with initial condition $u(\vv{x},0)=u_0(\vv{x})$ which is assumed to be compactly 
supported in the computational domain, $\Omega_i$.
Consider a rectangular cuboidal computational domain with boundary faces parallel to
one of the axes (see Fig.~\ref{fig:cuboidal}) with $\Omega_i =(x_l,x_r)\times\field{D}$
where $\field{D}=[-d_2,d_2)\times[-d_3,d_3)$, also termed as \emph{interior} domain. 
By using the decomposition of 
$u(\vv{x},t)\in \fs{L}^2(\field{R}\times\field{D})=\fs{L}^2(\Omega_i)
\oplus \fs{L}^2(\Omega_e)$ where $\Omega_{e}=(\field{R}\times\field{D})\setminus\overline{\Omega}_i$ 
also termed as \emph{exterior} domain. With this in mind, the equivalent formulation of 
problem~\eqref{eq:3D-SE} can be stated as  
\begin{equation}\label{eq:split-IVP-3D}
\begin{split}
\text{interior}:\;
&\left\{\begin{aligned}
&i\partial_tv+\triangle v=0,\;(\vv{x},t)\in\Omega_i\times\field{R}_+,\\
&v(\vv{x},0)= u_0(\vv{x}), \quad(\vv{x},t)\in\Omega_i,\\
\end{aligned}\right.\\
\text{exterior}:\;
&\left\{\begin{aligned}
&i\partial_tw+\triangle w=0,\;(\vv{x},t)\in\Omega_e\times\field{R}_+,\\
&w(\vv{x},0)=0,\quad\vv{x}\in\Omega_e,\\
&\lim_{|{x_1}|\rightarrow\infty} |w(\vv{x},t)| =0
\end{aligned}\right.\\
\text{continuity}:\; 
&v|_{\field{F}}=w|_{\field{F}},\; 
\partial_{n}v|_{\field{F}}=\partial_{n}w|_{\field{F}},
\end{split}
\end{equation}
where $\triangle=(\partial_{x_1}^2+\beta\partial_{x_2}^2+\beta\partial_{x_3}^2)$ 
and $\field{F}$ denotes the common regions of the interior and the exterior
domains. Also, the initial data is assumed to be compactly supported in $\Omega_i$.
Let us introduce the notations
\begin{equation}
\vv{x}_{\perp}=
\begin{pmatrix}
x_2\\
x_3
\end{pmatrix}
\in\field{D},\quad
\vs{\zeta}_{\perp}=
\begin{pmatrix}
\zeta_2,\\
\zeta_3
\end{pmatrix}
\end{equation}
where $\vs{\zeta}_{\perp}$ denote the covariable corresponding to $\vv{x}_{\perp}$ 
to be used in the two dimensional Fourier transform. Assuming the periodicity behaviour 
along $\vv{e}_2$ and $\vv{e}_3$, we can define the Fourier harmonics as
\begin{equation}
\field{H}_{\perp}
=\left\{\left.\pi\left(\frac{m_2}{d_2},\frac{m_3}{d_3}\right)^{\tp}\right|\;
m_2,m_3\in\field{Z}\right\}.
\end{equation}
We can express the field $w(x_1,x_2,x_3,t)$ in terms of a Fourier series as
\begin{equation}
\begin{split}
&w(x_1,x_2,x_3,t)=\sum_{\vs{\zeta}_{\perp}\in\field{H}_{\perp}}
\wtilde{w}_{\vs{\zeta}_{\perp}}(x_1,t)e^{i\vs{\zeta}_{\perp}\cdot\vv{x}_{\perp}},\\
&W(x_1,x_2,x_3,z)=\sum_{\vs{\zeta}_{\perp}\in\field{H}_{\perp}}
\wtilde{W}_{\vs{\zeta}_{\perp}}(x_1,z)e^{i\vs{\zeta}_{\perp}\cdot\vv{x}_{\perp}},
\end{split}
\end{equation}
where 
$\wtilde{W}_{\vs{\zeta}_{\perp}}(x_1,z)=\OP{L}_t[\wtilde{w}_{\vs{\zeta}_{\perp}}(x_1,t)]$ 
denotes the Laplace transform of the field.
Employing the Fourier expansion with respect to $x_2$ and $x_3$
and taking the Laplace transform of the right exterior problem yields
\begin{equation}
(\partial_{x_1}^2+\alpha^2)\wtilde{W}_{\vs{\zeta}_{\perp}}(x_1,z)=0,\quad x_1\in(x_r,\infty),
\end{equation}
where $\alpha=\sqrt{iz-\beta\|\vs{\zeta}_{\perp}\|^2}=\sqrt{iz-\beta{\zeta}_{\perp}^2}$ 
such that $\sqrt{\cdot}$ denotes the branch
with $\Im(\alpha)>0$. The boundedness of the solution leads to the relation
\begin{equation}
\partial_{x_1}\wtilde{W}_{\vs{\zeta}_{\perp}}(x_1,z) =i\alpha\wtilde{W}(x_1,z).
\end{equation}
To facilitate the inverse Laplace transform, we rewrite the above equation as
\begin{equation}\label{eq:dtn-map-freq-3d}
\OP{L}^{-1}[\partial_{x_1}\wtilde{W}_{\vs{\zeta}_{\perp}}]
=\OP{L}^{-1}[i\alpha^{-1}]\star\OP{L}^{-1}[\alpha^2\wtilde{W}_{\vs{\zeta}_{\perp}}],
\end{equation}
where `$\star$' represents the convolution operation.
Using the shifting property of Laplace transforms and introducing fractional operators 
of order $1/2$, we obtain the following result
\begin{equation}\label{eq:dtn-map-3d}
\begin{split}
\partial_{x_1}\wtilde{w}_{\vs{\zeta}_{\perp}}(x_1,t) 
&=e^{i\pi/4}e^{-i \beta\zeta_{\perp}^2t}\partial_t^{-1/2}e^{i\beta\zeta_{\perp}^2t}\\
&\quad \left[i\partial_t\wtilde{w}_{\vs{\zeta}_{\perp}}(x_1,t)
 -\beta{\zeta}_{\perp}^2\wtilde{w}_{\vs{\zeta}_{\perp}}(x_1,t)\right].
\end{split}
\end{equation}
In order to obtain the inverse Fourier transform of the above relation, let us 
introduce the kernel which is defined in a distributional sense as
\begin{equation}\label{eq:kernel-periodic-3d}
\mathcal{H}(\vv{x}_{\perp},t)=\frac{1}{4d_2d_3} \sum_{\vs{\zeta}_{\perp}\in\field{H}_{\perp}} 
e^{i(\vs{\zeta}_{\perp}\cdot\vv{x}_{\perp}-\beta\zeta_{\perp}^2t)}.
\end{equation}
Next by introducing the integral definition for the operator $\partial_t^{-1/2}$ and 
taking the inverse Fourier transform, the DtN map reads as
\begin{widetext}
\begin{equation}\label{eq:IPS-3D-TBC}
\begin{split}
\partial_{x_1}{w}(\vv{x},t)
&=\frac{e^{i\pi/4}}{\sqrt{\pi}}\int_0^t\int\int_{\field{R}^2}
\left[i\partial_{\tau}+\beta\triangle'_{\perp}\right]{w}(x_1,\vv{x}'_{\perp},\tau)
\frac{\mathcal{H}(\vv{x}_{\perp}-\vv{x}'_{\perp},t-\tau)}{\sqrt{t-\tau}}d^2\vv{x}'_{\perp}d\tau\\
&=-(\partial_{t}-i\beta\triangle_{\perp})\frac{e^{-i\pi/4}}{\sqrt{\pi}}
\int_0^t\int\int_{\field{R}^2}w(x_1,\vv{x}'_{\perp},\tau)
\frac{\mathcal{H}(\vv{x}_{\perp}-\vv{x}'_{\perp},t-\tau)}{\sqrt{t-\tau}}d^2\vv{x}'_{\perp}d\tau,
\end{split}
\end{equation}
where the last step is obtained using the following result from
the calculus of fractional operators~\cite{MR1993}:
\begin{equation}
\partial_t\left[e^{-i\beta\zeta_{\perp}^2t}\partial_t^{-1/2}e^{i\beta\zeta_{\perp}^2t}
\wtilde{w}_{\vs{\zeta}_{\perp}}(x_1,t)\right] =e^{-i\beta\zeta_{\perp}^2t}\partial_t^{-1/2}
e^{i\beta\zeta_{\perp}^2 t}\partial_t \wtilde{w}_{\vs{\zeta}_{\perp}}(x_1,t).
\end{equation}
In order to express the DtN map compactly, we introduce the notation 
\begin{equation}
(\partial_t-i\beta\triangle_{\perp})^{-1/2}f(\vv{x}_{\perp},t)\\
=\frac{1}{\sqrt{\pi}}\int_0^t\iint_{\field{R}^2}f(\vv{x}'_{\perp},\tau)
\frac{\mathcal{H}(\vv{x}_{\perp}-\vv{x}'_{\perp},t-\tau)}{\sqrt{t-\tau}}d^2\vv{x}'_{\perp}d\tau,
\end{equation}
\end{widetext}
so that
\[ (\partial_t-i\beta\triangle_{\perp})^{1/2}f
=(\partial_t-i\beta\triangle_{\perp})[(\partial_t-i\beta\triangle_{\perp})^{-1/2}f].\]
Using the continuity relations defined in~\eqref{eq:split-IVP}, the DtN
map on the right segment can be compactly written as
\begin{equation}
\partial_{x_1}{v}(\vv{x},t)
+e^{-i\pi/4}(\partial_t-i\beta\triangle_{\perp})^{1/2}v(\vv{x},t)=0.
\end{equation}
Set $\vv{d}_{\perp}=\left(d_2,d_3\right)^{\tp}$. An equivalent formulation of 
the IVP~\eqref{eq:2D-SE} on the computational
domain $\Omega_i = (x_l,x_r)\times\field{D}$
with periodic boundary condition along the axis $\vv{e}_2$ and $\vv{e}_3$ is given by
\begin{equation}\label{eq:3D-SE-CT}
\left\{\begin{aligned}
&i\partial_tu+(\partial_{x_1}^2+\beta\partial_{x_2}^2+\beta\partial_{x_3}^2)u=0,\; (\vv{x},t)\in\Omega_i\times\field{R}_+,\\
&u(\vv{x},0)=u_0(\vv{x})\in \fs{L}^2(\Omega_i),\quad\supp\, u_0\subset\Omega_i,\\
&\partial_{n}{u}+
e^{-i\pi/4}(\partial_{t}-i\beta\triangle_{\perp})^{1/2}u=0,
\; \vv{x}\in\field{F}_l\cup\field{F}_r,\\
&u(x_1,\vv{x}_{\perp}+2\vv{d}_{\perp},t)=u(x_1,\vv{x}_{\perp},t),\quad t>0,
\end{aligned}\right.
\end{equation}
where $\field{F}_l$ and $\field{F}_r$ denote the left and right boundary faces 
of the cuboidal domain, respectively.
\subsection{Fractional operators based approach}\label{sec:CT-CQ-3D}
\begin{figure}[!htbp]
\begin{center}
\def\myscale{1}
\includegraphics[scale=0.7]{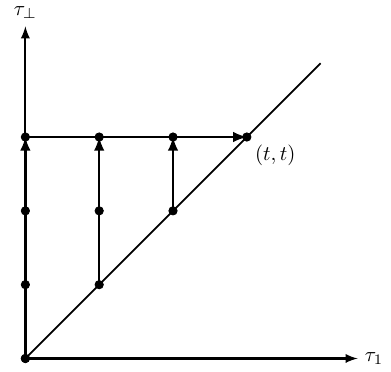}
\end{center}
\caption{\label{fig:IVP-auxi-periodic-cb} A schematic is shown in this figure which
depicts the evolution of the auxiliary field $\varphi(x_1,\vv{x}_{\perp},\tau_1,\tau_{\perp})$ 
in the $(\tau_1,\tau_{\perp})$-plane starting from the diagonal which also serves as 
initial conditions for solving the IVPs (arrow in the line depicts the 
evolution direction in time). Note that TBCs on $\partial{\Omega_i}$ require 
the history of the auxiliary field along the horizontal line up to the diagonal in the 
$(\tau_1,\tau_{\perp})$-plane.}
\end{figure}

In this section, we explore the possibility of writing the DtN maps in terms of 
(time-)fractional operators. Following the approach presented in 
Sec.~\ref{sec:CT-CQ}, let us recall the form of the DtN map given 
by~\eqref{eq:dtn-map-3d} in terms of the Fourier variables so that
\begin{equation}
\begin{split}
&\partial_{x_1}\wtilde{w}_{\vs{\zeta}_{\perp}}(x_1,t) \\
&=e^{i\pi/4}e^{-i \beta\zeta_{\perp}^2t}\partial_t^{-1/2}
 e^{i\beta\zeta_{\perp}^2t}\left[i\partial_t
 -\beta\zeta_{\perp}^2\right]\wtilde{w}_{\vs{\zeta}_{\perp}}(x_1,t)\\
&=ie^{i\pi/4}e^{-i \beta\zeta_{\perp}^2t}\partial_t^{1/2}
e^{i\beta\zeta_{\perp}^2t}\wtilde{w}_{\vs{\zeta}_{\perp}}(x_1,t)
\end{split}
\end{equation}
Taking inverse Fourier transform yields
\begin{equation*}
\begin{split}
&\partial_{x_1}{w}(\vv{x},t)\\
&=-{e^{-i\frac{\pi}{4}}} \sum_{\vs{\zeta}_{\perp}\in\field{H}_{\perp}}
e^{-i \beta\zeta_{\perp}^2t}\partial_t^{1/2}\left[e^{i\beta\zeta_{\perp}^2t}
\wtilde{w}_{\vs{\zeta}_{\perp}}(x_1,t)e^{i\vs{\zeta}_{\perp}\cdot\vv{x}_{\perp}} \right]\\
&=-{e^{-i\frac{\pi}{4}}}\partial^{1/2}_{t'}\sum_{\vs{\zeta}_{\perp}\in\field{H}_{\perp}}
\left[e^{-i\beta\zeta_{\perp}^2(t-t')}
\wtilde{w}_{\vs{\zeta}_{\perp}}(x_1,t')e^{i\vs{\zeta}_{\perp}\cdot\vv{x}_{\perp}}
\right]_{t'=t}.
\end{split}
\end{equation*}
Introducing the auxiliary function $\varphi(x_1,\vv{x}_{\perp},\tau_1,\tau_{\perp})$ such that
\begin{equation}
\begin{split}
&\varphi(x_1,\vv{x}_{\perp},\tau_1,\tau_{\perp})\\
&=\sum_{\vs{\zeta}_{\perp}\in\field{H}_{\perp}}
 e^{-i\beta\zeta_{\perp}^2(\tau_{\perp}-\tau_1)}\wtilde{w}_{\vs{\zeta}_{\perp}}(x_1,\tau_1)
 e^{i\vs{\zeta}_{\perp}\cdot\vv{x}_{\perp}}\\
&=\sum_{\vs{\zeta}_{\perp}\in\field{H}_{\perp}}
\left[ e^{+i\beta\zeta_{\perp}^2\tau_1}\wtilde{w}_{\vs{\zeta}_{\perp}}(x_1,\tau_1)\right]
e^{i(\vs{\zeta}_{\perp}\cdot\vv{x}_{\perp}-\beta\zeta_{\perp}^2 \tau_{\perp})}
\end{split}
\end{equation}
The DtN map on the boundary face $\field{F}_r$ can now be expressed compactly as
\begin{equation}
\begin{split}
&\partial_{x_1}{u}(\vv{x},t)=\left.-{e^{-i\frac{\pi}{4}}}
\partial^{1/2}_{\tau_1}\varphi(x_1,\vv{x}_{\perp},\tau_1,\tau_{\perp})\right|_{\tau_1,\tau_{\perp}=t}.
\end{split}
\end{equation}
By definition numerical implementation of this fractional derivative requires the history of 
the auxiliary function from the start of the computation which can be facilitated by 
the IVP satisfied by the auxiliary function: 
\begin{equation}\label{eq:FSE-auxi-ivp-3d}
\left\{\begin{aligned}
&(\tau_{\perp},\vv{x}_{\perp})\in(\tau_1,t]\times\field{D}:\\
&\quad [i\partial_{\tau_{\perp}}+\beta\triangle_{\perp}]
\varphi(x_1,\vv{x}_{\perp},\tau_1,\tau_{\perp})=0;\\
&\quad \varphi(x_1,\vv{x}_{\perp}+2\vv{d}_{\perp},\tau_1,\tau_{\perp})
=\varphi(x_1,\vv{x}_{\perp},\tau_1,\tau_{\perp}),\\
&\quad \varphi(x_1,\vv{x}_{\perp},\tau_1,\tau_1)=w(x_1,\vv{x}_{\perp},\tau_1).
\end{aligned}\right.
\end{equation}
Note that the history is needed for all $\tau_1\in[0,t]$ which means that we must 
solve the IVP for all $\tau_1\in[0,t]$ over $\field{F}_r$. It is interesting to note 
that the IVP satisfied by the auxiliary function on the face $\field{F}_r$ turns 
out to be a two-dimensional Schr\"{o}dinger equation. 
The computation of required history can be understood with the help of the
Fig.~\ref{fig:IVP-auxi-periodic-cb} where we note that the 
history of the field is needed along the horizontal line up to the diagonal in the
 $(\tau_1,\tau_{\perp})$-plane. The arrows point towards the direction of evolution
of the auxiliary function starting from the diagonal values which also serve as the 
initial conditions. 

The DtN maps on the boundary faces $\field{F}_l\cup\field{F}_r$ present 
in~\eqref{eq:3D-SE-CT} in terms of a time-fractional 
derivative operator reads as
\begin{equation}\label{eq:dtn-maps-cq-3d}
\partial_{n}{u} = - e^{-i\frac{\pi}{4}}\left.\partial^{1/2}_{\tau_1}
\varphi(x_1,\vv{x}_{\perp},\tau_1,\tau_{\perp})\right|_{\tau_1=\tau_{\perp}=t}.
\end{equation}
Here, we would like to mention that fractional operator based approach for the
numerical implementation of DtN maps is amenable to convolution quadratures but
this approach requires the history from the start of the computation that grows 
with each time-step making it expensive computationally as well as from a 
storage point of view.


\subsection{Novel Pad\'e based approach}\label{sec:CT-NP-3D}
\begin{figure}[!htbp]
\begin{center}
\includegraphics[scale=0.7]{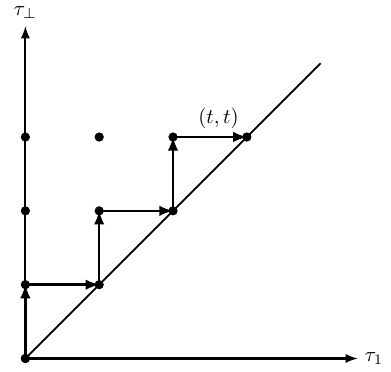}
\end{center}
\caption{\label{fig:IVP-auxi-pade-cb}A schematic depiction of the evolution of the 
auxiliary fields $\varphi_{k,a_1}(\vv{x}_{\perp},\tau_1,\tau_{\perp})$ in the 
$(\tau_1,\tau_{\perp})$-plane is provided in this figure. }
\end{figure}

In this section, we would like to explore the Pad\'e approximant based representation 
for the $1/2$-order temporal derivative in~\eqref{eq:dtn-maps-cq-3d} and its
other ramifications. This would allow us to obtain an effectively local approximation 
for the fractional derivative operator. We contrast this with the situation where 
the history grows with each time-step making it expensive computationally as well 
as from a storage point of view. 
Let us start by introducing a compact notation for the auxiliary functions as:
\begin{equation}
\varphi_{a_1}(\vv{x}_{\perp},\tau_1,\tau_{\perp})=\varphi(x_{a_1},\vv{x}_{\perp},\tau_1,\tau_{\perp})
,\quad a_1\in\{l,r\}.
\end{equation}
Let $z\in\field{C}$, then the Laplace transform of the auxiliary function can be denoted as  
\begin{equation}
\OP{L}_{\tau_1}[\varphi_{a_1}(\vv{x}_{\perp},\tau_1,\tau_{\perp})]
=\varPhi_{a_1}(\vv{x}_{\perp},z,\tau_{\perp}).
\end{equation}
The DtN maps described in~\eqref{eq:dtn-maps-cq-3d} on the faces $\field{F}_{a_1}$  
can be restated as 
\begin{equation}
\partial_{x_1}{u}(\vv{x},t)\pm e^{-i\pi/4}\left.\OP{L}^{-1}_{z}\left[
\sqrt{z}\varPhi_{a_1}\right]\right|_{\tau_1=\tau_{\perp}=t}=0.
\end{equation}
Inserting the rational approximation of order $M$ for the function $\sqrt{z}$ as
\begin{equation}\label{eq:maps-npade1-3d}
\partial_{x_1}{u}(\vv{x},t)\pm e^{-i\frac{\pi}{4}}\left.\OP{L}^{-1}_{z}\left[
R_M(z)\varPhi_{a_1}\right]\right|_{\tau_1=\tau_{\perp}=t}\approx 0.
\end{equation}
The map expressed above is still non-local, however, the rational nature of the 
expression allows us to obtain an effectively local approximation for the DtN 
operation by introducing the auxiliary fields $\varphi_{k,a_1}(\vv{x}_{\perp},\tau_1,\tau_{\perp})$ 
on the boundary faces $\field{F}_{a_1}$ as
\begin{equation}
(z+\eta^2_k)^{-1}\varPhi_{a_1}(\vv{x}_{\perp},z,\tau_{\perp})=\varPhi_{k,a_1}(\vv{x}_{\perp},z,\tau_{\perp}),
\end{equation}
where $k=1,2,\ldots,M$. In the physical space, every $\varphi_{k,a_1}$ satisfy the IVP given by
\begin{equation}\label{eq:auxi-npade-3d}
(\partial_{\tau_1}+\eta^2_k)\varphi_{k,a_1}(\vv{x}_{\perp},\tau_1,\tau_{\perp})
=\varphi(\vv{x}_{\perp},\tau_1,\tau_{\perp}),
\end{equation}
with the initial conditions set to $\varphi_{k,a_1}(\vv{x}_{\perp},0,\tau_{\perp})=0$. The DtN map 
now reads as 
\begin{equation}\label{eq:maps-npade-3d}
\partial_{n_1}{u}+e^{-i\frac{\pi}{4}}\left[ b_0 u(\vv{x},t) -
\sum_{k=1}^M b_k\varphi_{k,a_1}(\vv{x}_{\perp},t,t)\right]=0.
\end{equation}
The solution to the ODEs~\eqref{eq:auxi-npade-3d} work out to be
\begin{equation*}
\varphi_{k,a_1}(\vv{x}_{\perp},\tau_1,\tau_{\perp})
=\int_0^{\tau_1}e^{-\eta_k^2(\tau_1-s)}\varphi(x_1,\vv{x}_{\perp},s,\tau_{\perp})ds.
\end{equation*}
It is easy to verify that the auxiliary fields $\phi_{k,a_1}(\vv{x}_{\perp},\tau_1,\tau_{\perp})$ 
satisfy the IVPs given by
\begin{equation}\label{eq:ivp-pade-auxi-3d}
[i\partial_{\tau_{\perp}}+\beta\triangle_{\perp}]\varphi_{k,a_1}(\vv{x}_{\perp},\tau_1,\tau_{\perp})=0.
\end{equation}
The evolution of the auxiliary fields $\varphi_{k,a_1}(\vv{x}_{\perp},\tau_1,\tau_{\perp})$
in the $(\tau_1,\tau_2)$-plane can be understood with the help of the schematic 
shown in Fig.~\ref{fig:IVP-auxi-pade-cb}. The diagonal to diagonal computation 
of the fields $\varphi_{k,a_1}(\vv{x}_{\perp},\tau_1,\tau_{\perp})$ consists of first advancing 
the fields using the IVP established in~\eqref{eq:ivp-pade-auxi-3d} and then using 
the ODEs in~\eqref{eq:auxi-npade-3d} for the second movement. The novel Pad\'e 
approach makes the numerical scheme efficient from a storage point of view 
which is also clear from the schematic presented. 

\section{Numerical Implementation--2D}\label{sec:numerical-implementation-2d}
In this section, we address the complete numerical solution 
of the initial boundary-value problem (IBVP) stated in~\eqref{eq:2D-SE-CT}.
The formulation of the discrete linear system for the IBVP requires 
temporal as well as spatial discretization of the problem. 
To this end, we first discretize the temporal derivative using one-step 
methods, namely, backward differentiation formula of order 1 (BDF1) and the 
trapezoidal rule (TR). Subsequently, a compatible temporal discretization 
scheme is developed for the boundary conditions where we strive to formulate 
the novel boundary maps at the time-discrete level as a Robin-type 
boundary condition for each time-step. Finally, for the spatial discretization
of the resulting problem, we use a Legendre-Galerkin method where we construct
a boundary-adapted basis in terms of the Legendre polynomials followed by a 
boundary lifting process in order to enforce the boundary conditions exactly.

For the computational domain $\Omega_i$, we introduce a reference domain 
$\Omega_i^{\text{ref.}}=\field{I}\times\pi\field{I}$ keeping in view the typical 
domain of definition of the orthogonal polynomials being used.
In order to describe the associated linear maps between the reference domain and 
the actual computational domain, we introduce the variables 
$y_1,y_2\in\Omega_i^{\text{ref.}}$ such that
\begin{equation*}
\begin{split}
&x_1 = J_1 y_1+\bar{x}_1,\quad J_1 = \frac{1}{2}(x_r-x_l),\quad
\bar{x}_1=\frac{1}{2}(x_l+x_r), \\
&x_2 = J_2 y_2,\quad J_2 = \frac{d}{\pi},\quad \beta_1=J_1^{-2},\quad
\beta_2=\beta J_2^{-2}.
\end{split}
\end{equation*}
Let $\Delta t$ denote the time-step. 
In the rest of the section, with a slight abuse of notation, we switch to the 
variables $y_1\in\field{I}$ and $y_2\in\pi\field{I}$ for all the discrete approximations 
of the dependent variables. For instance, $u^j(\vv{y})$ is taken to approximate 
$u(\vv{x},j\Delta t)$ for $j=0,1,2,\ldots,N_t-1$. The temporal discretization of the 
interior problem using one-step methods is enumerated below:
\begin{itemize}
\item The BDF1 based discretization is given by
\begin{equation}\label{eq:2d-se-bdf1}
\begin{split}
&i\frac{u^{j+1}-u^{j}}{\Delta t}+\triangle u^{j+1}=0, \quad \rho = 1/\Delta t,\\
&\implies\left(\beta_1\partial^2_{y_1}
+\beta_2\partial^2_{y_2}\right)u^{j+1}+i\rho u^{j+1}
=i\rho u^j.
\end{split}
\end{equation}
\item The TR based discretization is given by
\begin{equation}\label{eq:2d-disc-tr}
\begin{split}
&i\frac{u^{j+1}-u^{j}}{\Delta t}+\triangle u^{j+1/2}=0,\quad \rho = 2/\Delta t,\\
&\implies\left(\beta_1\partial^2_{y_1}+\beta_2\partial^2_{y_2}\right)v^{j+1}+i\rho v^{j+1}
=i\rho u^j,
\end{split}
\end{equation}
where we have used the staggered samples of the field defined as
\begin{equation}
v^{j+1} = u^{j+1/2}=\frac{u^{j+1}+u^j}{2},\quad v^0=0.
\end{equation}
\end{itemize}
The rest of this section is organized as follows: Sec.~\ref{sec:dtn-maps}
addresses the temporal discretization of the boundary maps which is followed 
by a discussion of the resulting spatial problem in Sec.~\ref{sec:IBVP}.
\subsection{Discretizing the boundary conditions}\label{sec:dtn-maps}
As discussed above, we aim at a compatible temporal discretization of the novel 
boundary conditions formulated in earlier sections in order to solve the IBVP. To 
this end, we first discuss the convolution quadrature (CQ) based discretization 
of the $1/2$-order time-fractional derivative present in the TBCs in 
Sec.~\ref{sec:CT-CQ}. Here, we try to match the time-stepping methods in the CQ 
based discretization of the boundary conditions to that of the interior problem 
(namely, BDF1 and TR) to preserve the overall stability of the numerical recipe~\cite{PM2015}.
Next, we discuss the temporal discretization of the boundary maps developed in 
Sec.~\ref{sec:CT-NP} and Sec.~\ref{sec:CT-CP} using Pad\'{e} approximants 
based representation for the fractional operator present in the TBCs. Lastly,
we present a CQ based discretization of the fractional derivative operator present 
in the boundary maps obtained under high-frequency approximation in Sec.~\ref{sec:CT-HF}.

In order to facilitate the temporal discretization of the TBCs, we introduce the 
equispaced samples of the auxiliary function, $\varphi^{j,k}(y_1,y_2)$, to 
approximate $\varphi(x_1,x_2,j\Delta t,k\Delta t)$ for $j,k=0,1,\ldots N_t-1$. It 
follows from the definition of the auxiliary function that the diagonal values of 
the auxiliary function at the discrete level are determined by the interior field,
i.e., $\varphi^{j,j}(y_1,y_2) = u^j(y_1,y_2)$.

Let $\Delta t$ denote the time-step so that the temporal grid becomes 
$t_j=j\Delta t,\;j\in\field{N}_0$. The fractional operators 
$\partial_t^{\pm1/2}$ can be implemented using the 
LMM based CQ method. For a given LMM, let the CQ weights be denoted by 
$\{\omega^{(\pm1/2)}_{k}|\;k\in\field{N}_0\}$, then
\begin{equation}
[\partial_t^{\pm1/2}F_{\pm}]^{j+1}
=\rho^{\pm 1/2}F_{\pm}^{j+1}
+\rho^{\pm 1/2}\sum_{k=1}^{j}\omega^{(\pm1/2)}_{j+1-k}F_{\pm}^k,
\end{equation}
where $\rho$ is a function of $\Delta t$ and it is chosen such that 
$\omega^{(\pm1/2)}_{0}=1$. In the following we restrict ourselves to the 
trapezoidal rule (TR) and backward differentiation formulae (BDF) of order $1$. 
The recipe for computing the CQ weights for each of the time-stepping methods 
is enumerated below:
\begin{itemize}
\item CQ--BDF1 : The discretization scheme for the time-fractional operators is said to be 
`CQ-BDF1' if the underlying time-stepping method used in designing the quadrature 
is BDF1. Let $\nu\in\{ +1/2, -1/2\} $ and $\rho =1/\Delta t$, then the quadrature weights 
can be computed as follows 
\begin{equation}
\omega^{(\nu)}_k=\left[{\left(k-1-\nu \right)}/{k}\right]\omega^{(\nu)}_{k-1},\quad k\geq 1,
\end{equation}
where $\omega^{(\nu)}_0=1$.
\item CQ--TR : Let $\rho =2/\Delta t$, then the quadrature weights for the TR based
CQ scheme can be computed as
\begin{equation}
\begin{split}
&\omega^{(1/2)}_k=(-1)^j\omega^{(-1/2)}_k, \\
&\omega^{(-1/2)}_k=
\begin{cases}
C_{k/2},      & k\;\text{even},\\
C_{(k-1)/2},  & k\;\text{odd},\\
\end{cases},\\
&\text{where}\quad C_n=\frac{1\cdot 3\cdots(2n-1)}{n!2^n}.
\end{split}
\end{equation}
The quadrature weights can also be
generated using the following recurrence relation~\cite{SV2023}:
\begin{equation}
(k+1)\omega^{(\nu)}_{k+1}=(k-1)\omega^{(\nu)}_{k-1}-2\nu\omega^{(\nu)}_k,\,k\geq2,
\end{equation}
where $\omega^{(\nu)}_0=1$ and $\omega^{(\nu)}_1=-2\nu$.
\end{itemize}

\subsubsection{CQ--BDF1}
The discrete scheme for the boundary maps is said to be `CQ--BDF1' if the 
underlying one-step method in the CQ scheme is BDF1. Let 
$(\omega^{(1/2)}_k)_{k\in\field{N}_0}$ denote the corresponding quadrature weights, 
then the resulting discretization for the DtN maps described 
in~\eqref{eq:dtn-maps-cq} in terms of $\varphi^{j,k}(y_1,y_2)$ stated as a 
Robin-type boundary condition reads as
\begin{multline}
\partial_{y_1}u^{j+1}\pm\sqrt{\rho/\beta_1}e^{-i\pi/4}u^{j+1}\\
=\mp\sqrt{\rho/\beta_1}e^{-i\pi/4}
\OP{B}\left[\{\varphi^{k,j+1}(y_{a_1},y_2)\}_{k=0}^j\right],
\end{multline}
where
\begin{multline}
\mathcal{B}^{j+1}_{a_1}(y_2)
=\OP{B}\left[\{\varphi^{k,j+1}(y_{a_1},y_2)\}_{k=0}^j\right]\\
=\sum_{k=1}^{j+1}\omega^{(1/2)}_{k}\varphi^{j+1-k,j+1}(y_{a_1},y_2),\quad
a_1\in\{r,l\}.
\end{multline}
Set $\alpha_k=\sqrt{\rho/\beta_k}e^{-i\pi/4}$ for $k=1,2$. The maps in the
compact form reads as
\begin{equation}\label{eq:cq-maps-bdf1}
\begin{split}
&\partial_{y_1}u^{j+1}-\alpha_1 u^{j+1}
=+\alpha_1\mathcal{B}^{j+1}_l(y_2),\\
&\partial_{y_1}u^{j+1}+\alpha_1 u^{j+1}
=-\alpha_1\mathcal{B}^{j+1}_r(y_2).
\end{split}
\end{equation}
To obtain the history of the field $\varphi(x_1,x_2,\tau_1,\tau_2)$ needed to compute the 
boundary sums, we employ Fourier Galerkin method (on account of periodicity along $x_2$)
to solve the IVP~\eqref{eq:FSE-auxi-ivp}. For the space discretization of the auxiliary equation, 
we use the Fourier-Galerkin method which employs the series
\begin{equation}
\varphi^{j,k}(y_1,y_2)
=\sum_{m\in\field{J}_2}\wtilde{\varphi}_m^{j,k}(y_1,y_2)e^{im y_2},
\end{equation}
where the index set is given by
\begin{equation}
\field{J}_2=\left\{-\frac{N}{2},\ldots,\frac{N}{2}-1\right\}.
\end{equation}
BDF1-based discretization for the IVP in~\eqref{eq:FSE-auxi-ivp} reads as follows: 
\begin{multline}
i\frac{\wtilde{\varphi}^{k,j+1}_m-\wtilde{\varphi}^{k,j}_m}{\Delta t}
-\beta_2m^2 \wtilde{\varphi}^{k,j+1}_m=0,\\
\implies\wtilde{\varphi}^{k,j+1}_m=\frac{\rho}{(\rho+i\beta_2m^2)}\wtilde{\varphi}^{k,j}_m,
\quad m\in\field{J}_2.
\end{multline}
In terms of $\alpha_2=\sqrt{\rho/\beta_2}e^{-i\pi/4}$, the coefficients $\wtilde{\varphi}^{k,j+1}_m$ 
satisfy the recurrence relation
\begin{equation}
\wtilde{\varphi}^{k,j+1}_m=\frac{\alpha_2^2}{(\alpha_2^2+m^2)}\wtilde{\varphi}^{k,j}_m,
\end{equation}
with $\wtilde{\varphi}^{j,j}_m=\wtilde{u}^{j}_m$.
\subsubsection{CQ--TR}
The discrete scheme is said to be `CQ--TR' if the underlying one-step method
in the CQ scheme is TR. Let $(\omega^{(1/2)}_k)_{k\in\field{N}_0}$ denote the
corresponding quadrature weights, then the resulting discretization scheme for 
the DtN map described in~\eqref{eq:dtn-maps-cq} in terms of $\varphi^{j,k}(y_1,y_2)$ 
stated as a Robin-type boundary condition (consistent
with staggered samples of interior field) reads as
\begin{multline}
\partial_{y_1}v^{j+1} \pm\alpha_1v^{j+1}\\
=\mp\frac{1}{2}\alpha_1
\left(\OP{B}\left[ \{\varphi^{k,j+1}_{a_1}\}_{k=0}^j\right]
+\OP{B}\left[\{\varphi^{k,j}_{a_1}\}_{k=0}^{j-1}\right]\right),
\end{multline}
where the history operator $\OP{B}$ is defined in the similar manner as in the
case of BDF1 except the quadrature weights.
Let $ \mathcal{B}^{j+1/2}=(\mathcal{B}^{j+1}+\mathcal{B}^{j})/2$, then the maps 
can be compactly written as
\begin{equation}\label{eq:cq-maps-tr}
\begin{split}
&\partial_{y_1}v^{j+1}-\alpha_1 v^{j+1}=+\alpha_1\mathcal{B}^{j+1/2}_l(y_2),\\
&\partial_{y_1}v^{j+1}+\alpha_1 v^{j+1}=-\alpha_1\mathcal{B}^{j+1/2}_r(y_2).
\end{split}
\end{equation}
Once again, we solve the IVP for the auxiliary function in similar manner as that of 
BDF1 except using TR-based discretization of the IVP~\eqref{eq:FSE-auxi-ivp} which
reads as follows: 
\begin{multline}
i\left(\frac{\wtilde{\varphi}^{k,j+1}_m-\wtilde{\varphi}^{k,j}_m}{\Delta t}\right)
-\beta_2m^2\left(\frac{\wtilde{\varphi}^{k,j+1}_m+\wtilde{\varphi}^{k,j}_m}{2}\right)=0,\\
\implies\wtilde{\varphi}^{k,j+1}_m=\left(\frac{\rho-i\beta_2m^2}{\rho+i\beta_2m^2}\right)
\wtilde{\varphi}^{k,j}_m.
\end{multline}
In terms of $\alpha_2=\sqrt{\rho/\beta_2}e^{-i\pi/4}$, the coefficients $\wtilde{\varphi}^{j,k}_m$ 
satisfy the recurrence relation
\begin{equation}
\wtilde{\varphi}^{k,j+1}_m=\left(\frac{\alpha_2^2-m^2}{\alpha_2^2+m^2}\right)
\wtilde{\varphi}^{k,j}_m.
\end{equation}
\subsubsection{NP-BDF1}
The discrete scheme for the boundary maps is said to be `NP--BDF1' if the 
temporal derivatives in the realization of the boundary maps, dubbed as `NP', 
developed in Sec.~\ref{sec:CT-NP} are discretized using the one-step method 
BDF1. To achieve this, we start by writing the time-discrete form of the DtN 
maps present in~\eqref{eq:maps-npade} in the reference domain as
\begin{equation}\label{eq:map-npade-bdf1}
\sqrt{\beta_1}\partial_{y_1}\wtilde{u}_m^{j+1}
\pm e^{-i\pi/4}\left[ b_0\wtilde{u}_m^{j+1} -
\sum_{k=1}^M b_k\wtilde{\varphi}^{j+1,j+1}_{k,m}\right]=0.
\end{equation}
In order to turn these equations into a Robin-type boundary condition, we need to 
compute the discrete samples $\varphi^{j+1,j+1}_{k,m}$ 
by considering the ODEs established for these auxiliary function earlier. 
Here, we discuss the BDF1-based discretization of~\eqref{eq:auxi-npade} by first
employing Fourier expansion: 
Setting $\rho=1/\Delta t$ and for $k=1,2,\ldots,M$, we have
\begin{multline}
\frac{\wtilde{\varphi}^{j+1,j+1}_{k,m}-\wtilde{\varphi}^{j, j+1}_{k,m}}{\Delta t} 
+\eta^2_k\wtilde{\varphi}^{j+1,j+1}_{k,m}
=\wtilde{\varphi}_m^{j+1,j+1}\\
\implies\wtilde{\varphi}^{j+1,j+1}_{k,m} 
=\frac{\rho}{(\rho+\eta^2_k)}\wtilde{\varphi}^{j,j+1}_{k,m}
+\frac{1}{(\rho+\eta^2_k)}\wtilde{u}_m^{j+1}.
\end{multline}
Introduce the scaling $\ovl{\eta}_k=\eta_k/\sqrt{\rho}$ so that
\begin{equation}
\wtilde{\varphi}^{j+1,j+1}_{k,m} 
=\frac{1}{(1+\ovl{\eta}^2_k)}\left[\wtilde{\varphi}^{j,j+1}_{k,m}
+\frac{1}{\rho}\wtilde{u}_m^{j+1}\right].
\end{equation}
The value of $\wtilde{\varphi}^{j,j+1}_{k,m}$ can be obtained by propagation 
along $\tau_2$ using $\wtilde{\varphi}^{j,j}_{k,m}$ as the initial condition:
\begin{equation}
\wtilde{\varphi}^{j,j+1}_{k,m}
=\frac{1}{(1+\alpha_2^{-2}m^2)}\wtilde{\varphi}^{j,j}_{k,m},
\quad m\in\field{J}_2.
\end{equation}
Next we plug-in the value of $\wtilde{\varphi}^{j+1,j+1}_{k,m}$ in~\eqref{eq:map-npade-bdf1}
to obtain
\begin{multline}
\partial_{y_1}\wtilde{u}_m^{j+1}\pm\alpha_1\varpi\wtilde{u}_m^{j+1}\\
=\mp\frac{\alpha_1}{(1+\alpha_2^{-2} m^2)}\left[\sum_{k=1}^M \Gamma_k\wtilde{\varphi}^{j,j}_{k,m}\right]
=\mp\alpha_1\wtilde{\mathcal{B}}^{j+1}_{m,a_1},
\end{multline}
where we have set $\ovl{b}_0=b_0/\sqrt{\rho},\; \ovl{b}_k =b_k/\sqrt{\rho}$ and
\begin{equation}
\varpi =\ovl{b}_0 +\frac{1}{\rho}\sum_{k=1}^M\Gamma_k,\quad
\Gamma_{k} =-{\ovl{b}_{k}}/{(1+\ovl{\eta}^2_{k})}.
\end{equation}
Note that these Robin type maps are written for each component in the Fourier series. 
Also observe that $\varpi$ is independent of $m$ in this case which allows us to take 
inverse Fourier transform of these maps and obtain the form of maps in the physical space 
as follows:
\begin{equation}\label{eq:npade-maps-bdf1}
\partial_{y_1} u^{j+1}\pm\alpha_1\varpi u^{j+1}=\mp\alpha_1\mathcal{B}^{j+1}_{a_1},\\
\end{equation}
where
\begin{equation}
\mathcal{B}^{j+1}_{a_1}(y_2)=\sum_{m\in\field{J}_2}\wtilde{\mathcal{B}}^{j+1}_{m,a_1}e^{imy_2}
,\quad a_1\in\{r,l\}.
\end{equation}
\subsubsection{NP--TR}
The discrete scheme for the boundary maps is said to be `NP--TR' if the 
temporal derivatives in the realization of the boundary maps, dubbed as `NP', 
developed in Sec.~\ref{sec:CT-NP} are discretized using the one-step method 
TR. To achieve this, we start by writing the time-discrete form of the DtN 
maps present in~\eqref{eq:maps-npade} in the reference domain as
\begin{equation}\label{eq:map-npade-tr}
\sqrt{\beta_1}\partial_{y_1}\wtilde{v}_m^{j+1}
\pm e^{-i\pi/4}\left[ b_0\wtilde{v}_m^{j+1} -
\sum_{k=1}^M b_k\wtilde{\varphi}^{j+1/2,j+1/2}_{k,m}\right]=0.
\end{equation}
In order to turn these equations into a Robin-type boundary condition, we need to 
compute the discrete samples $\varphi^{j+1/2,j+1/2}_{k,m}$ 
by considering the ODEs established for these auxiliary function earlier. 
Here, we discuss the TR-based discretization of~\eqref{eq:auxi-npade} by first
employing Fourier expansion: 
Setting $\rho=2/\Delta t$ and for $k=1,2,\ldots,M$, we have
\begin{multline}
\frac{\wtilde{\varphi}^{j+1,j+1}_{k,m}-\wtilde{\varphi}^{j,j+1}_{k,m}}{\Delta t} 
+\eta^2_k\frac{\wtilde{\varphi}^{j+1,j+1}_{k,m}+\wtilde{\varphi}^{j,j+1,j}_{k,m}}{2}\\
= \frac{\wtilde{\varphi}_m^{j+1,j+1}+\wtilde{\varphi}_m^{j,j+1}}{2}.
\end{multline}
Introduce the scaling $\ovl{\eta}_k=\eta_k/\sqrt{\rho}$ so that
\begin{multline}
\wtilde{\varphi}^{j+1,j+1}_{k,m} 
=\frac{(1-\ovl{\eta}^2_k)}{(1+\ovl{\eta}^2_k)}\wtilde{\varphi}^{j,j+1}_{k,m}\\
+\frac{2/\rho}{(1+\ovl{\eta}^2_k)}\left[\wtilde{v}_m^{j+1}
+\frac{\wtilde{\varphi}_m^{j,j+1}-\wtilde{\varphi}_m^{j,j}}{2}\right].
\end{multline}
For the staggered configuration, we have
\begin{multline}\label{eq:map-novel-pade-trap}
\wtilde{\varphi}^{j+1/2,j+1/2}_{k,m} 
=\frac{1}{2}\left[\left(\frac{1-\ovl{\eta}^2_k}{1+\ovl{\eta}^2_k}\right)
\wtilde{\varphi}^{j,j+1}_{k,m}+\wtilde{\varphi}^{j,j}_{k,m}\right]\\
+\frac{1/\rho}{(1+\ovl{\eta}^2_k)}\left[\wtilde{v}_m^{j+1}
+\frac{\wtilde{\varphi}_m^{j,j+1}-\wtilde{\varphi}_m^{j,j}}{2}\right].
\end{multline}
The values of $\wtilde{\varphi}^{j,j+1}_{k,m}$ and $\wtilde{\varphi}^{j,j+1}_{m}$ 
can be obtained by propagation along $\tau_2$ using
$\wtilde{\varphi}^{j,j}_{k,m}$ and $\wtilde{\varphi}^{j,j}_{m}$ as 
the initial condition:
\begin{multline}
\wtilde{\varphi}^{j,j+1}_{k,m}
=\frac{\left(1-\alpha_2^{-2}m^2\right)}{\left(1+\alpha_2^{-2}m^2\right)}
\wtilde{\varphi}^{j,j}_{k,m},\\
\frac{1}{2}\left[\wtilde{\varphi}^{j,j+1}_{m}-\wtilde{\varphi}^{j,j}_{m}\right]
=\frac{-\alpha_2^{-2}m^2}{\left(1+\alpha_2^{-2}m^2\right)}
\wtilde{\varphi}^{j,j}_{m},\quad m\in\field{J}_2.
\end{multline}
Next we plug-in the value of $\wtilde{\varphi}^{j+1/2,j+1/2}_{k,m}$
in~\eqref{eq:map-npade-tr} which leads to following compact form of the 
boundary maps
\begin{equation}
\partial_{y_1}\wtilde{v}_m^{j+1}\pm\alpha_1\varpi\wtilde{v}_m^{j+1}
=\mp\alpha_1\wtilde{\mathcal{B}}^{j+1/2}_{m,a_1},\\
\end{equation}
where $\varpi$ is same as that in BDF1 and the history functions are given by
\begin{equation}
\begin{split}
&\wtilde{\mathcal{B}}^{j+1/2}_{m,a_1}\\
&=\sum_{k=1}^M\frac{-\ovl{b}_k}{2}\left[\left(\frac{1-\ovl{\eta}^2_k}{1+\ovl{\eta}^2_k}\right)
\left(\frac{1-\alpha_2^{-2} m^2}{1+\alpha_2^{-2}m^2}\right)+1\right]\wtilde{\varphi}^{j,j}_{k,m}\\
&\quad+\sum_{k=1}^M\frac{\Gamma_k}{\rho}\left(\frac{-\alpha_2^{-2}m^2}{1+\alpha_2^{-2}m^2}\right)
\wtilde{u}^{j}_{m}.
\end{split}
\end{equation}
Note that these Robin type maps are written for each component in Fourier series. 
Also $\varpi$ is independent of $m$ in this case as well which allows us to take 
inverse Fourier transform of these maps as:
\begin{equation}
\partial_{y_1} v^{j+1}\pm\alpha_1\varpi v^{j+1}=\mp\alpha_1\mathcal{B}^{j+1/2}_{a_1},\\
\end{equation}
where
\begin{equation}
\mathcal{B}^{j+1/2}_{a_1}(y_2)=\sum_{J_2}\wtilde{\mathcal{B}}^{j+1/2}_{m,a_1}e^{imy_2}
,\quad a_1\in\{r,l\}.
\end{equation}
\subsubsection{CP-BDF1}
The discrete scheme for the boundary maps is said to be `CP--BDF1' if the 
temporal derivatives in the realization of the boundary maps, dubbed as `CP', 
developed in Sec.~\ref{sec:CT-CP} are discretized using the one-step method 
BDF1. To achieve this, we start by writing the time-discrete form of the DtN 
maps present in~\eqref{eq:maps-cpade} in the reference domain as
\begin{equation}\label{eq:map-cpade-bdf1}
\sqrt{\beta_1}\partial_{y_1}\wtilde{u}_m^{j+1}
\pm e^{-i\pi/4}\left[ b_0\wtilde{u}_m^{j+1} -
\sum_{k=1}^M b_k\wtilde{\varphi}^{j+1}_{k,m}\right]=0.
\end{equation}
In order to turn these equations into a Robin-type boundary condition, we need to 
compute the discrete samples $\varphi^{j+1}_{k,m}$ by considering the ODEs established
for these auxiliary function earlier. 
Here, we discuss the BDF1-based discretization of~\eqref{eq:auxi-cpade} by first
employing Fourier expansion: 
Setting $\rho=1/\Delta t$ and for $k=1,2,\ldots,M$, we have
\begin{multline}
i\frac{\wtilde{\varphi}^{j+1}_{k,m}-\wtilde{\varphi}^{j}_{k,m}}{\Delta t} 
-(\beta_2m^2-i\eta^2_k)\wtilde{\varphi}^{j+1}_{k,m}
= i\wtilde{u}_m^{j+1}\\
\wtilde{\varphi}^{j+1}_{k,m} 
=\frac{1}{(\rho+i\beta_2 m^2+\eta^2_k)}
\left[\rho\wtilde{\varphi}^j_{k,m}+\wtilde{u}_m^{j+1}\right].
\end{multline}
Introduce the scaling $\ovl{\eta}_k=\eta_k/\sqrt{\rho}$ so that
\begin{equation}
\wtilde{\varphi}^{j+1}_{k,m}
=\frac{1}{1+\alpha_2^{-2}m^2+\ovl{\eta}_k^2}\left[\wtilde{\varphi}^j_{k,m}
+\frac{1}{\rho}\wtilde{u}_m^{j+1}\right],
\end{equation}
where $\wtilde{\varphi}^{j}_{k,m}$ and $\wtilde{u}^j_m$ approximate the Fourier
coefficients $\OP{F}_{y_2}[\varphi_k(y_2,j\Delta t)]$ and
$\OP{F}_{y_2}[u(y_1,y_2,j\Delta t)]$, respectively. 
Next we plug-in the value of $\wtilde{\varphi}^{j+1}_{k,m}$ in~\eqref{eq:map-cpade-bdf1} 
to obtain the following boundary maps
\begin{multline}\label{eq:pade-maps-bdf1}
\partial_{y_1}\wtilde{u}_m^{j+1}\pm \alpha_1\varpi_{m}\wtilde{u}_m^{j+1}\\
= \mp\alpha_1\sum_{k=1}^M\Gamma_{k,m}\wtilde{\varphi}^{j}_{k,m} 
=\mp\alpha_1\wtilde{\mathcal{B}}_{m,a_1}^{j+1},
\end{multline}
where we have set $\ovl{b}_0=b_0/\sqrt{\rho},\; \ovl{b}_k =b_k/\sqrt{\rho}$ and
\begin{equation}
\begin{split}
\varpi_{m}   &= \ovl{b}_0 + \frac{1}{\rho}\sum_{k=1}^M\Gamma_{k,m},\\
\Gamma_{k,m} &=\frac{-\ovl{b}_k}{(1+\alpha_2^{-2}m^2+\ovl{\eta}_k^2)}.
\end{split}
\end{equation}
Note that these Robin type maps are written for each component in Fourier series. 
Also observe that $\varpi$ is $m$-dependent in this case which poses a difficulty in 
taking inverse Fourier transform of these maps to express them in physical space. 
\subsubsection{CP-TR}
The discrete scheme for the boundary maps is said to be `CP--TR' if the 
temporal derivatives in the realization of the boundary maps, dubbed as `CP', 
developed in Sec.~\ref{sec:CT-CP} are discretized using the one-step method 
TR. To achieve this, we start by writing the time-discrete form of the DtN 
maps present in~\eqref{eq:maps-cpade} in the reference domain as
\begin{equation}\label{eq:map-cpade-tr}
\sqrt{\beta_1}\partial_{y_1}\wtilde{v}_m^{j+1}
\pm e^{-i\pi/4}\left[ b_0\wtilde{v}_m^{j+1} -
\sum_{k=1}^M b_k\wtilde{\varphi}^{j+1/2}_{k,m}\right]=0.
\end{equation}
In order to turn these equations into a Robin-type boundary condition, we need to 
compute the discrete samples $\varphi^{j+1/2}_{k,m}$ by considering the ODEs established
for these auxiliary function earlier. 
Here, we discuss the TR-based discretization of~\eqref{eq:auxi-cpade} by first
employing Fourier expansion: 
Setting $\rho=2/\Delta t$ and for $k=1,2,\ldots,M$, we have
\begin{multline*}
i\frac{\wtilde{\varphi}^{j+1}_{k,m}-\wtilde{\varphi}^j_{k,m}}{\Delta t} 
-(\beta_2m^2-i\eta^2_k)\frac{\wtilde{\varphi}^{j+1}_{k,m}+\wtilde{\varphi}^j_{k,m}}{2}
=i\wtilde{v}_m^{j+1}\\
\wtilde{\varphi}^{j+1}_{k,m} 
=\frac{(\rho-i\beta_2m^2-\eta^2_k)}{(\rho+i\beta_2m^2+\eta^2_k)}\wtilde{\varphi}^j_{k,m}
+\frac{2\wtilde{v}_m^{j+1}}{(\rho+i\beta_2m^2+\eta^2_k)}.
\end{multline*}
Introduce the scaling $\ovl{\eta}_k=\eta_k/\sqrt{\rho}$ so that
\begin{equation*}
\wtilde{\varphi}^{j+1}_{k,m} 
=\frac{(1-\alpha_2^{-2}m^2-\ovl{\eta}^2_k)}{(1+\alpha_2^{-2}m^2+\ovl{\eta}^2_k)}
\wtilde{\varphi}^j_{k,m}+\frac{(2/\rho)\wtilde{v}_m^{j+1}}{(1+\alpha_2^{-2}m^2+\ovl{\eta}^2_k)}.
\end{equation*}
For the staggered configuration, we have
\begin{equation*}\label{eq:map-pade-trap}
\wtilde{\varphi}^{j+1/2}_{k,m} 
= \frac{1}{(1+\alpha_2^{-2}m^2+\ovl{\eta}_k^2)}
\left[\wtilde{\varphi}^j_{k,m}
+(1/\rho)\wtilde{v}_m^{j+1}\right].
\end{equation*}
Next we plug-in the value of $\wtilde{\varphi}^{j+1/2}_{k,m}$ in~\eqref{eq:map-cpade-tr}
to obtain the following boundary maps
\begin{equation}\label{eq:disc-map-pade}
\partial_{y_1}\wtilde{v}_m^{j+1}\pm\alpha_1\varpi_m\wtilde{v}_m^{j+1}
=\mp\alpha_1\wtilde{\mathcal{B}}^{j+1/2}_{m,a_1},\\
\end{equation}
where $\varpi_m$ is same as CP--BDF1 and the history function is given by
\begin{equation}
\wtilde{\mathcal{B}}^{j+1/2}_{m,a_1}
=\left[\sum_{k=1}^M\Gamma_{k,m}\wtilde{\varphi}^{j}_{k,m}\right]_{y_1=y_{a_1}},
\quad a_1\in\{l,r\}.
\end{equation}
Also observe that $\varpi$ is $m$-dependent in this case which poses a difficulty in 
taking inverse Fourier transform of these maps to express them in 
physical space. 
\subsubsection{HF-BDF1}
The discrete scheme for the boundary maps obtained under high-frequency approximation
is said to be `HF--CQ--BDF1' if the underlying one-step method in the CQ scheme is BDF1.
Let $\nu\in\{-1/2,+1/2\}$ and $(\omega^{(\nu)}_k)_{k\in\field{N}_0}$ denote the
corresponding quadrature weights, then the resulting discretization of the DtN maps 
described in~\eqref{eq:dtn-maps-hf} 
reads as
\begin{multline}
\sqrt{\beta_1}\partial_{y_1}u^{j+1}\pm\sqrt{\rho}e^{-i\pi/4}u^{j+1}
\mp\frac{e^{i\pi/4}}{2\sqrt{\rho}}\beta_2\partial^2_{y_2}u^{j+1}\\
=\mp\sqrt{\rho}e^{-i\pi/4}\mathcal{B}^{j+1}_{a_1,1/2}
\pm\frac{e^{i\pi/4}}{2\sqrt{\rho}}\beta_2\partial^2_{y_2}\mathcal{B}^{j+1}_{a_1,-1/2},
\end{multline}
where history functions are given by
\[
\mathcal{B}^{j+1}_{a_1,\nu} =\sum_{k=1}^{j+1}\omega^{(\nu)}_{k}u^{j+1-k},
\quad a_1\in\{l,r\}.
\]
Set $\alpha_j=\sqrt{\rho/\beta_j}e^{-i\pi/4}$ for $j=1,2$, then the DtN maps
can be expressed compactly as
\begin{multline}
\partial_{y_1}u^{j+1}\pm\alpha_1 u^{j+1}\mp\frac{1}{2}\alpha_1\alpha_2^{-2}\partial^2_{y_2}u^{j+1}
\\
=\mp\alpha_1\mathcal{B}^{j+1}_{a_1,1/2}
\pm\frac{1}{2}\alpha_1\alpha_2^{-2}\partial^2_{y_2}\mathcal{B}^{j+1}_{a_1,-1/2}.
\end{multline}
On account of periodicity along $y_2$, we expand the field as a Fourier series 
and rewrite the DtN maps for each component as
\begin{multline}
\partial_{y_1}\wtilde{u}_m^{j+1}\pm\alpha_1 \wtilde{u}_m^{j+1}
\mp\frac{1}{2}\alpha_1\alpha_2^{-2}m^2\wtilde{u}_m^{j+1}\\
=\mp\alpha_1\wtilde{\mathcal{B}}^{j+1}_{m,a_1,1/2}
\pm\frac{1}{2}\alpha_1\alpha_2^{-2}m^2\wtilde{\mathcal{B}}^{j+1}_{m,a_1,-1/2}.
\end{multline}
Under high-frequency approximation, we can express DtN maps as a Robin-type map as 
\begin{equation}\label{eq:maps-bdf1-hf}
\begin{split}
&\partial_{y_1}\wtilde{u}_m^{j+1}-\alpha_1\varpi_m\wtilde{u}_m^{j+1}
=+\alpha_1\wtilde{\mathcal{B}}^{j+1}_{m,l},\\
&\partial_{y_1}\wtilde{u}_m^{j+1}+\alpha_1\varpi_m\wtilde{u}_m^{j+1}
=-\alpha_1\wtilde{\mathcal{B}}^{j+1}_{m,r},
\end{split}
\end{equation}
where we introduced $\varpi_m=1+\alpha_2^{-2}m^2/2 $ and
\begin{equation}
\wtilde{\mathcal{B}}^{j+1}_{m,a_1}
=\wtilde{\mathcal{B}}^{j+1}_{m,a_1,1/2}+\frac{1}{2}\alpha_2^{-2}m^2
\wtilde{\mathcal{B}}^{j+1}_{m,a_1,-1/2}.
\end{equation}

\subsubsection{HF-TR}
The discrete scheme for the boundary maps obtained under high-frequency approximation
is said to be `HF--CQ--TR' if the underlying one-step method in the CQ scheme is TR. 
Let $\nu\in\{-1/2,+1/2\}$ and  $(\omega^{(\nu)}_k)_{k\in\field{N}_0}$ denote the
corresponding quadrature weights, then the resulting discretization of the DtN maps 
(consistent with staggered samples
of the field) described in~\eqref{eq:dtn-maps-hf} reads as
\begin{multline}
\sqrt{\beta_1}\partial_{y_1}v^{j+1}\pm\sqrt{\rho}e^{-i\pi/4}v^{j+1}
\mp\frac{e^{i\pi/4}}{2\sqrt{\rho}}\beta_2\partial^2_{y_2}v^{j+1}\\
=\mp\sqrt{\rho}e^{-i\pi/4}\mathcal{B}^{j+1/2}_{a_1,1/2}
\pm\frac{e^{i\pi/4}}{2\sqrt{\rho}}\beta_2\partial^2_{y_2}\mathcal{B}^{j+1/2}_{a_1,-1/2},
\end{multline}
where the history functions are given by
\begin{equation}
\mathcal{B}^{j+1/2}_{a_1,\nu}
=\sum_{k=1}^{j+1}\omega^{(\nu)}_{k}v^{j+1-k},\quad a_1\in\{l,r\}.
\end{equation}
Following along the similar lines as discussed in BDF1 case, we can obtain a compact
Robin-type boundary map as
\begin{equation}\label{eq:maps-tr-hf}
\begin{split}
&\partial_{y_1}\wtilde{v}_m^{j+1}-\alpha_1\varpi_m\wtilde{v}_m^{j+1}
=+\alpha_1\wtilde{\mathcal{B}}^{j+1/2}_{m,l},\\
&\partial_{y_1}\wtilde{v}_m^{j+1}+\alpha_1\varpi_m\wtilde{v}_m^{j+1}
=-\alpha_1\wtilde{\mathcal{B}}^{j+1/2}_{m,r},
\end{split}
\end{equation}
where $\varpi_m=1+\alpha_2^{-2}m^2/2$ and for $a_1\in\{r,l\}$ 
\begin{equation}
\wtilde{\mathcal{B}}^{j+1/2}_{m,a_1}
=\wtilde{\mathcal{B}}^{j+1/2}_{m,a_1,1/2}+\frac{1}{2}\alpha_2^{-2}m^2
 \wtilde{\mathcal{B}}^{j+1/2}_{m,a_1,-1/2}.
\end{equation}
Note that $\varpi_m$ is $m$-dependent which makes it difficult to take inverse 
Fourier transform of the maps obtained in \eqref{eq:maps-tr-hf} in order to 
express them in physical space. This is similar to the case of conventional-Pad\'e
approach. 
\subsection{Numerical Solution of the IBVP}\label{sec:IBVP}
\subsubsection{CQ and NP}
Having discussed the temporal discretization of the novel boundary maps considered 
in this work together with that of the interior problem, we would like to
address the spatial problem at hand. For spatial discretization, we use a method which is 
a combination of the Legendre and the Fourier basis within the framework of Galerkin 
formulation. We construct a boundary-adapted basis in terms of Legendre polynomials 
which makes it convenient to enforce the boundary conditions and, at the same time, 
yields a banded linear system. Note that the formulation of the linear system requires
a homogeneous form of boundary maps which is achieved by a boundary lifting process. 
Set $\field{I}=[-1,1]$ and let $L_n(y)$ denote the Legendre polynomial of degree $n$. 
Define the polynomial space 
\begin{equation}
\fs{P}_N = \Span\left\{L_p(y)|\;p=0,1,\ldots,N,\;y\in\field{I}\right\}.
\end{equation}
Similarly, using the Fourier basis, we define the linear vector space 
\begin{equation}
\fs{V}_N=\Span\left\{\psi_k=\exp(iky)|\;k\in\field{J}_2,\;y\in\pi\field{I}\right\}.
\end{equation}
where the index set is given by
\begin{equation}
\field{J}_2=\left\{-\frac{N}{2},\ldots,\frac{N}{2}-1\right\}.
\end{equation}
For the 2D problem, we consider the tensor product space given by 
$\fs{P}_N\otimes\fs{V}_N$. In order to enforce the 
boundary conditions exactly, we introduce the space of boundary-adapted basis as follows:
\begin{equation*}
\fs{X}_N=\left\{u\in\fs{P}_N\otimes\fs{V}_N\left|\;
\begin{aligned}
&(\partial_{y_1}-\kappa)u(y_1,y_2)|_{y_1=-1}=0,\\
&(\partial_{y_1}+\kappa)u(y_1,y_2)|_{y_1=+1}=0
\end{aligned}\right.
\right\},
\end{equation*}
where $\kappa$ is determined by the Robin-type formulation of the TBCs specific to the method
chosen for the temporal discretization of the boundary maps.
\begin{rem}
The form of $\kappa$ follows from the Robin-type formulation of the discretized
version of the DtN-maps as discussed in Sec.~\ref{sec:dtn-maps}.
For CQ methods $\kappa=\alpha_1=\sqrt{\rho/\beta_1}\exp(-i\pi/4)$ while 
for the NP methods, we have $\kappa=\alpha_1\varpi$.
\end{rem}

Note that the CQ and NP methods afford the possibility to be casted in this form
on account of the fact that $\varpi$ is independent of the $m$-th Fourier component.
But CP and HF methods cannot be casted in this form as they have $m$-dependent
$\varpi$. The construction of the boundary-adapted basis and boundary lifting process 
will be discussed for these methods separately. Let the index set  $\{0,1,\ldots,N-2\}$ 
be denoted by $\field{J}_1$. For CQ and NP methods, we introduce a function $\chi$ to 
convert the discrete TBCs to a homogeneous form by stating the original field as
$u^j=w^j+\chi^{j}(y_1,y_2)$ where $w^j\in\fs{X}_N$ so that
\begin{equation}\label{eq:1D-SE-homogenized}
\begin{split}
[-\alpha_1^{-2}\partial^2_{y_1}-\alpha_2^{-2}\partial^2_{y_2}+1](w^{j+1}+\chi^{j+1})
= u^j,\\
\begin{aligned}
&\left.(\partial_{y_1}-\kappa)w^{j+1}\right|_{y_1=-1}=0,\\
&\left.(\partial_{y_1}+\kappa)w^{j+1}\right|_{y_1=+1}=0,
\end{aligned}
\end{split}
\end{equation}
and the field $\chi^j(y)$ is forced to satisfy the constraints below
\begin{equation}
\begin{split}
&\left.(\partial_{y_1}-\kappa)\chi^{j}(y_1,y_2)\right|_{y_1=-1}=+\alpha_1\mathcal{B}^j_l(y_2),\\
&\left.(\partial_{y_1}+\kappa)\chi^{j}(y_1,y_2)\right|_{y_1=+1}=-\alpha_1\mathcal{B}^j_r(y_2).
\end{split}
\end{equation}
The form of the equations suggest the use of an ansatz of the form
$\chi^j(y_1,y_2)=c^j_l\chi^j_l+c^j_r\chi^j_r$ with
$c^j_l=\alpha_1\mathcal{B}^j_l$ and $c^j_r=-\alpha_1\mathcal{B}^j_r$ so that the
lifting functions satisfy
\begin{equation}
\begin{split}
&\left.(\partial_{y_1}-\kappa)
\begin{pmatrix}
\chi_l\\
\chi_r
\end{pmatrix}
\right|_{y_1=-1}=
\begin{pmatrix}
1\\
0
\end{pmatrix},\\
&\left.(\partial_{y_1}+\kappa)
\begin{pmatrix}
\chi_l\\
\chi_r
\end{pmatrix}
\right|_{y_1=+1}=
\begin{pmatrix}
0\\
1
\end{pmatrix}.
\end{split}
\end{equation}
From here, the degree of the lifting functions $\chi_{l,r}$ can be
inferred to be utmost $1$. Exploiting the fact 
that we can expand them in terms of Legendre polynomials, we have
\begin{equation}
\begin{pmatrix}
	\chi_l(y_1)\\
	\chi_r(y_1)
\end{pmatrix}=A_0
\begin{pmatrix}
	L_0(y_1)\\
	L_1(y_1)
\end{pmatrix},
\end{equation}
where $A_0$ is the unknown transformation matrix. Solving the linear system for $A_0$ yields
\begin{equation}
\left\{\begin{aligned}
\chi_l(y_1) &=-\frac{1}{2\kappa}L_0(y_1)+\frac{1}{2(\kappa+1)}L_1(y_1),\\
\chi_r(y_1) &=+\frac{1}{2\kappa}L_0(y_1)+\frac{1}{2(\kappa+1)}L_1(y_1).
\end{aligned}\right.
\end{equation}
The lifted field then takes the form
\begin{equation}\label{eq:lifted-field}
u^j=w^j
+\alpha_1\left[\chi_l(y_1)\mathcal{B}^j_{l}(y_2)-\chi_r(y_1)\mathcal{B}^j_{r}(y_2)\right],
\end{equation}
which can be expanded as
\begin{multline}
u^j=w^j-\frac{\alpha_1}{2\kappa}\sum_{q\in\field{J}_2}
\left[\wtilde{\mathcal{B}}^j_{q,r}+\wtilde{\mathcal{B}}^j_{q,l}\right]L_0(y_1)\psi_q(y_2)\\
-\frac{\alpha_1}{2(1+\kappa)}\sum_{q\in\field{J}_2}
\left[\wtilde{\mathcal{B}}^j_{q,r}-\wtilde{\mathcal{B}}^j_{q,l}\right]L_1(y_1)\psi_q(y_2),
\end{multline}
where we have employed the Fourier expansion of the history functions
\begin{equation}
\mathcal{B}^j_{a_1}(y_2)=\sum_{q\in\field{J}_2}\wtilde{\mathcal{B}}^j_{q,a_1}\psi_q(y_2),
\quad a_1\in\{l,r\}.
\end{equation}
Next, we choose the basis $\theta_{p,q}\in\fs{X}_N$ defined as
\begin{equation}
\theta_{p,q}(y_1,y_2)=\phi_p(y_1)\psi_q(y_2),\quad p\in\field{J}_1,\;q\in\field{J}_2,
\end{equation}
where $\{\phi_{p}|\,p\in\field{J}_1\}$ can be constructed using the ansatz 
$\phi_{p}(y_1)=L_p(y_1)+a_pL_{p+1}(y_1)+b_pL_{p+2}(y_1)$.
Imposing the boundary conditions in the definition of $\fs{X}_N$, the 
sequences $(a_p)_{p\in\field{J}_1}$ and $(b_p)_{p\in\field{J}_1}$ work out to be
\begin{equation}
a_p=0,\quad
b_p=-\frac{\kappa+\frac{1}{2}p(p+1)}{\kappa+\frac{1}{2}(p+2)(p+3)},
\quad p\in\field{J}_1.
\end{equation}
The testing functions for the variational formulation in context of a Galerkin method
can now be stated as
\begin{equation}
\wtilde{\theta}_{p,q}(y_1,y_2)=\phi_{p}(y_1)\psi^*_{q}(y_2)
=\phi_{p}(y_1)\psi_{-q}(y_2).
\end{equation}
The mass matrix and the stiffness matrix for the boundary-adapted Legendre basis are
denoted by $M_1 = (\melem{m}^{(1)}_{k,j})_{k,j\in\field{J}_1}$ and 
$S_1 = (\melem{s}^{(1)}_{k,j})_{k,j\in\field{J}_1}$, respectively, where
\begin{equation}\label{eq:sys-mat-2d}
\begin{split}
&\melem{s}^{(1)}_{j,k}=-(\phi_j,\phi''_k)_{\field{I}}
=\begin{cases}
-2(2k+3)b_k,&j=k,\\
0,&\mbox{otherwise},
\end{cases}\\
&\melem{m}^{(1)}_{j,k}=(\phi_j,\phi_k)_{\field{I}}
=\begin{cases}
\frac{2b_{k-2}}{2k+1},&j=k-2,\\
\frac{2}{2k+1}+\frac{2b_k^2}{2k+5},&j=k,\\
\frac{2b_{k}}{2k+5},&j=k+2,\\
0&\mbox{otherwise}.
\end{cases}
\end{split}
\end{equation}
The mass matrix and the stiffness matrix for the Fourier basis are denoted by 
$M_2 = (\melem{m}^{(2)}_{k,j})_{k,j\in\field{J}_2}$ and 
$S_2 = (\melem{s}^{(2)}_{k,j})_{k,j\in\field{J}_2}$, respectively, where
\begin{equation}
\begin{split}
&\melem{s}^{(2)}_{j,k}=-(\psi''_j,\psi^{*}_k)_{\pi\field{I}}
=\begin{cases}
2\pi k^2&,k=j\\
0&,\mbox{otherwise},
\end{cases}\\
&\melem{m}^{(2)}_{j,k}=(\psi_j,\psi^{*}_k)_{\pi\field{I}}
=\begin{cases}
2\pi,& k=j,\\
0,&\mbox{otherwise}.
\end{cases}
\end{split}
\end{equation}
Within the variational formulation, our goal is to find the 
approximate solution $u_N\in\fs{X}_N$ to~\eqref{eq:2d-se-bdf1} such that
\begin{multline}\label{eq:vform-SE2D}
-\alpha_1^{-2}\left(\partial^2_{y_1}u_N^{j+1},\wtilde{\theta}_{p,q}\right)_{\Omega_i^{\text{ref}}}
-\alpha^{-2}_2\left(\partial^2_{y_2}u_N^{j+1},\wtilde{\theta}_{p,q}\right)_{\Omega_i^{\text{ref}}}\\
+\left(u_N^{j+1},\wtilde{\theta}_{p,q}\right)_{\Omega_i^{\text{ref}}}
=\left(u_N^j,\wtilde{\theta}_{p,q}\right)_{\Omega_i^{\text{ref}}},
\end{multline}
for every $\wtilde{\theta}_{p,q}\in\fs{X}_{N}$ 
where $(u,\wtilde{\theta})_{\Omega_i^{\text{ref}}}=\int_{\Omega_i^{\text{ref}}}u\wtilde{\theta}d^2\vv{y}$
is the scalar product in $\fs{L}^2(\Omega_i^{\text{ref}})$. We approximate $u$ by polynomial 
interpolation over the Legendre-Gauss-Lobatto (LGL) nodes in the context of discrete 
Legendre transforms. 
The variational formulation~\eqref{eq:vform-SE2D} in terms of the 
unknown field $w^{j+1}$ (dropping the subscript `$N$' for the sake of brevity) 
can now be stated as
\begin{multline}\label{eq:v-form}
-\left[\alpha_1^{-2}\left(\partial^2_{y_1}w^{j+1},\wtilde{\theta}_{p,q}\right)
+\alpha_2^{-2}\left(\partial^2_{y_2}w^{j+1},\wtilde{\theta}_{p,q}\right)\right]\\
+\left(w^{j+1},\wtilde{\theta}_{p,q}\right)
=\left(u^j,\wtilde{\theta}_{p,q}\right)
+\mathcal{I},
\end{multline}
where
\begin{multline}
\mathcal{I}
=-\left(\chi^{j+1},\wtilde{\theta}_{p,q}\right)
+\underbrace{\alpha_1^{-2}\left(\partial^2_{y_1}\chi^{j+1},\wtilde{\theta}_{p,q}\right)}_{=0}\\
+\alpha_2^{-2}\left(\partial^2_{y_2}\chi^{j+1},\wtilde{\theta}_{p,q}\right).
\end{multline}
The non-zero term above works out to be
\begin{multline}
\partial^2_{y_2}\chi^{j+1}=
-\alpha_1\chi_r(y_1)\partial^2_{y_2}{\mathcal{B}}^{j+1}_r(y_2)\\
+\alpha_1\chi_l(y_1)\partial^2_{y_2}{\mathcal{B}}^{j+1}_l(y_2).
\end{multline}
Integration by parts for $a_1\in\{l,r\}$ yields
\begin{equation}
\left(\chi_{a_1}\partial^2_{y_2}\mathcal{B}^{j+1}_{a_1},\wtilde{\theta}_{p,q}\right)
=-q^2\left(\chi_{a_1}\mathcal{B}^{j+1}_{a_1},\wtilde{\theta}_{p,q}\right),
\end{equation}
so that
\begin{equation}
\mathcal{I}
=\alpha_1(1+\alpha_2^{-2}q^2)\left[(\mathcal{B}^{j+1}_{r}\chi_r,\wtilde{\theta}_{p,q})
-(\mathcal{B}^{j+1}_{l}\chi_l,\wtilde{\theta}_{p,q})\right].
\end{equation}
Let $\wtilde{u}_{p,q}$ denote the expansion coefficients of the field in Legendre-Fourier basis 
and $\what{w}_{p,q}$ denote the expansion coefficients in the boundary-adapted basis as
\begin{equation}
\begin{split}
w^{j}(\vv{y})
&=\sum_{p\in\field{J}_1,\;q\in\field{J}_2}\what{w}^{j}_{p,q}\theta_{p,q}(\vv{y}),\\
&=\sum_{p=0}^N\sum_{q\in\field{J}_2}\wtilde{w}^{j}_{p,q}L_{p}(y_1)\psi_{q}(y_2),\\
u^{j}(\vv{y})
&=\sum_{p=0}^N\sum_{q\in\field{J}_2}
\wtilde{u}^{j}_{p,q}L_{p}(y_1)\psi_{q}(y_2).
\end{split}
\end{equation}
Let us introduce the matrices 
$\what{W}=\left(\what{w}_{p,q}\right)$
and $\wtilde{U}=\left(\wtilde{u}_{p,q}\right)$ 
for convenience. In matrix form, the expansion coefficients for the 
fields in the Legendre basis assume the form
\begin{equation}
\begin{split}
\wtilde{U}^j=\wtilde{W}^j
&-\frac{\alpha_1}{2\kappa}\vv{e}_0\otimes
\left[\wtilde{\vs{\mathcal{B}}}^j_{r}+\wtilde{\vs{\mathcal{B}}}^j_{l}\right]^{\tp}\\
&-\frac{\alpha_1}{2(1+\kappa)}\vv{e}_1\otimes
\left[\wtilde{\vs{\mathcal{B}}}^j_{r}-\wtilde{\vs{\mathcal{B}}}^j_{l}\right]^{\tp},
\end{split}
\end{equation}
where $\wtilde{\vs{\mathcal{B}}}^j_{a_1},\;a_1\in\{l,r\},$ are the corresponding 
Fourier coefficients. The inner-products in the variational 
form involving the stiffness-matrix follows from
\begin{equation}
\begin{split}
&-\left[\alpha_1^{-2}\left(\partial^2_{y_1}w^{j+1},\wtilde{\theta}_{p,q}\right)
+\alpha_2^{-2}\left(\partial^2_{y_2}w^{j+1},\wtilde{\theta}_{p,q}\right)\right]\\
&=\sum_{p',q'}\left[
\alpha_1^{-2}\melem{s}^{(1)}_{p,p'}\what{w}^{j+1}_{p',q'}\melem{m}^{(2)}_{q',q}
+\alpha_2^{-2}\melem{m}^{(1)}_{p,p'}\what{w}^{j+1}_{p',q'}\melem{s}^{(2)}_{q',q}
\right]\\
&=\left(\alpha_1^{-2}S_1\what{W}^{j+1}M_2+\alpha_2^{-2}M_1\what{W}^{j+1}S_2\right)_{p,q},
\end{split}
\end{equation}
and that involving only the mass-matrices follow from
\begin{equation}
\begin{split}
(w^{j+1},\wtilde{\theta}_{p,q})
&=\sum_{p',q'}\melem{m}^{(1)}_{p,p'}\what{w}^{j+1}_{p',q'}\melem{m}^{(2)}_{q',q}\\
&=\left(M_1\what{W}^{j+1}M_2\right)_{p,q}.
\end{split}
\end{equation}
We can further simplify the expression for the linear system by collecting the 
right hand side of~\eqref{eq:v-form} into one term. To this end, let $f^j(y_1,y_2)$ 
be such that 
$(u^j,\wtilde{\theta}_{p,q})+\mathcal{I}=(f^j,\wtilde{\theta}_{p,q})$,
then it follows that 
\begin{multline}
\frac{\alpha_1}{4\pi}\sum_{q\in\field{J}_2}(1+\alpha_2^{-2}q^2)
\left(\frac{g^j_{q,+}}{\kappa}+\frac{g^j_{q,-}}{1+\kappa}\right)L_0(y_1)\psi_q(y_2)\\
+u^j=f^j,
\end{multline}
where $g^j_{q,\pm}=
\left(\mathcal{B}^{j+1}_{r}\pm\mathcal{B}^{j+1}_{l},\psi^*_{q}\right)_{\pi\field{I}}$.
Introducing the diagonal matrix  
$\mathcal{D}_2=\diag\left(\left\{1+\alpha_2^{-2}q^2\right\}_{q\in\field{J}_2}\right)$,
the discrete representation of $f^j$ in the matrix form works out to be
\begin{equation}
\begin{split}
\wtilde{F}^j=\wtilde{U}^j
&+\frac{\alpha_1}{2\kappa}\vv{e}_0\otimes
\left[\wtilde{\vs{\mathcal{B}}}^{j+1}_{r}+\wtilde{\vs{\mathcal{B}}}^{j+1}_{l}\right]^{\tp}\mathcal{D}_2\\
&+\frac{\alpha_1}{2(1+\kappa)}\vv{e}_1\otimes
\left[\wtilde{\vs{\mathcal{B}}}^{j+1}_{r}-\wtilde{\vs{\mathcal{B}}}^{j+1}_{l}\right]^{\tp}\mathcal{D}_2.
\end{split}
\end{equation}
The inner product $(f^j,\wtilde{\theta}_{p,q})$ can be computed from the
knowledge of the Legendre coefficients of the field $\wtilde{F}^j$ using the
specific form of the boundary-adapted basis. This can be achieved via
\emph{quadrature matrix} given in~\ref{app:ip}.
The linear system for~\eqref{eq:v-form} thus becomes
\begin{multline}\label{eq:linear-2d-system}
\alpha_1^{-2}S_1\what{W}^{j+1}M_2+\alpha^{-2}_2M_1\what{W}^{j+1}S_2\\
+M_1\what{W}^{j+1}M_2
=Q_1\Gamma\wtilde{F}^{j}M_2,
\end{multline}
which is recasted as follows
\begin{equation*}
\left(\alpha_1^{-2}M_2\otimes S_1+\alpha_2^{-2}S_2\otimes M_1
+M_2\otimes M_1\right)\what{\vv{w}}^{j+1}=\what{\vv{f}}^j, 
\end{equation*}
where $\otimes$ denotes the tensor product and $\what{\vv{w}},\;\what{\vv{f}}$
denote the column vectors obtained by stacking columns of matrices 
$\what{W},\;\what{F}$, respectively, below one another. The mass and stiffness 
matrices are defined in~\eqref{eq:sys-mat-2d} and the linear system is then 
solved using LU-decomposition method. 

\subsubsection{CP and HF}
For CP and HF methods, we have the following form of the TBCs casted as 
Robin-type boundary maps:
\begin{equation}
\partial_{y_1}\wtilde{u}_m^{j+1}\pm \alpha_1\varpi_m\wtilde{u}_m^{j+1}
=\mp\alpha_1\wtilde{\mathcal{B}}^{j+1}_{m,a_1},\quad a_1\in\{r,l\}.
\end{equation}
Note that these maps have $m$-dependent $\varpi$ which makes it difficult to take the 
inverse Fourier transform in order to express them in physical space. On account of
this fact, we will have to discuss the boundary lifting and basis construction for 
each Fourier component separately. In order to enforce these $m$-dependent boundary 
maps exactly, we define the space of boundary-adapted basis as
\begin{equation*}
\fs{W}_N=\left\{\wtilde{u}_m\in\fs{P}_N\left|\;
\begin{aligned}
&(\partial_{y_1}-\kappa_m)\wtilde{u}_m(y_1)|_{y_1=-1}=0,\\
&(\partial_{y_1}+\kappa_m)\wtilde{u}_m(y_1)|_{y_1=+1}=0
\end{aligned}\right.
\right\},
\end{equation*}
where $\kappa_m$ will be determined by the Robin-type formulation of TBCs specific to
CP and HF methods.
We introduce a function $\chi_m$ to convert the discrete BCs to a homogeneous form by
stating the original field as $\wtilde{u}_m^j=\wtilde{w}_m^j+\chi_m^{j}(y_1)$
where $\wtilde{w}_m^j\in\fs{W}_N$ so that
\begin{equation}
\begin{split}
[-\alpha_1^{-2}\partial^2_{y_1}+\alpha_2^{-2}m^2
+1](\wtilde{w}_m^{j+1}+\chi_m^{j+1}) = \wtilde{u}_m^j,\\
\begin{aligned}
&\left.(\partial_{y_1}-\kappa_m)\wtilde{w}_m^{j+1}\right|_{y_1=-1}=0,\\
&\left.(\partial_{y_1}+\kappa_m)\wtilde{w}_m^{j+1}\right|_{y_1=+1}=0,
\end{aligned}
\end{split}
\end{equation}
and the field $\chi_m^j(y_1)$ is forced to satisfy the constraints:
\begin{equation}
\begin{split}
&(\partial_{y_1}-\kappa_m)\chi_m^{j}(y_1)|_{y_1=-1}=+\alpha_1\wtilde{\mathcal{B}}^j_{m,l},\\
&(\partial_{y_1}+\kappa_m)\chi_m^{j}(y_1)|_{y_1=+1}=-\alpha_1\wtilde{\mathcal{B}}^j_{m,r},
\end{split}
\end{equation}
Following in the similar manner as in CQ and NP case, the lifting functions for
this case works out to be
\begin{equation}
\left\{\begin{aligned}
\chi_{m,l}(y_1) &=-\frac{1}{2\kappa_m}L_0(y_1)+\frac{1}{2(\kappa_m+1)}L_1(y_1),\\
\chi_{m,r}(y_1) &=+\frac{1}{2\kappa_m}L_0(y_1)+\frac{1}{2(\kappa_m+1)}L_1(y_1).
\end{aligned}\right.
\end{equation}
The lifted field becomes
\begin{equation}\label{eq:lifted-field-m}
\wtilde{u}_m^j({y_1})=\wtilde{w}_m^j({y_1})
+\alpha_1\chi_{m,l}\wtilde{\mathcal{B}}^j_{m,l}-\alpha_1\chi_{m,r}\wtilde{\mathcal{B}}^j_{m,r},
\end{equation}
which can be expanded as
\begin{multline}
\wtilde{u}_m^j(y_1)=\wtilde{w}_m^j(y_1)-\frac{\alpha_1}{2\kappa_m}
\left[\wtilde{\mathcal{B}}^j_{q,r}+\wtilde{\mathcal{B}}^j_{q,l}\right]L_0(y_1)\\
-\frac{\alpha_1}{2(1+\kappa_m)}
\left[\wtilde{\mathcal{B}}^j_{q,r}-\wtilde{\mathcal{B}}^j_{q,l}\right]L_1(y_1).
\end{multline}
Next we address the construction of the basis $\theta_{p,m}\in\fs{W}_N$ for the boundary-adapted 
space $\fs{W}_N$ as
\begin{equation}
\theta_{p,m}(y_1)=\phi_{p,m}(y_1),\quad p\in\field{J}_1,\;m\in\field{J}_2,
\end{equation}
where $\phi_{p,m}$ can be computed using the ansatz 
$\phi_{p,m}(y_1)=L_p(y_1)+a_{p,m}L_{p+1}(y_1)+b_{p,m}L_{p+2}(y_1)$.
Imposing the boundary conditions in the definition, the 
sequences $\{a_{p,m}\}$ and $\{b_{p,m}\}$ work out to be
\begin{equation}
a_{p,m}=0,\;
b_{p,m}=-\frac{\alpha_1\varpi_m+\frac{1}{2}p(p+1)}{\alpha_1\varpi_m+\frac{1}{2}(p+2)(p+3)}.
\end{equation}
The testing functions for the variational formulation in context of a Galerkin method
can now be stated as
\begin{equation}
\wtilde{\theta}_{p,m}(y_1)=\phi_{p,-m}(y_1)=\phi_{p,m}(y_1),
\end{equation}
where we used $\varpi_{-m}=\varpi_m$ and $p\in\field{J}_1,\;m\in\field{J}_2$.
The mass matrix and the stiffness matrix for the boundary adapted Legendre basis
(specific to each $m$-th Fourier component) are denoted by 
$\wtilde{M}_{m,1} = (\melem{m}^{(1)}_{p,p',m})_{p,p'\in\field{J}_1}$ and 
$\wtilde{S}_{m,1} = (\melem{s}^{(1)}_{p,p',m})_{p,p'\in\field{J}_1}$, respectively, where
\begin{equation}
\begin{split}
\melem{s}^{(1)}_{p,p',m}
&=-(\phi_{p',m},\phi''_{p,m})_{\field{I}}\\
&=\begin{cases}
-2(2p+3)b_{p,m},&p'=p,\\
0,&\mbox{otherwise}
\end{cases},\\
\melem{m}^{(1)}_{p,p',m}
&=(\phi_{p',m},\phi_{p,m})_{\field{I}}\\
&=\begin{cases}
\frac{2b_{p-2,m}}{2p+1},&p'=p-2,\\
\frac{2}{2p+1}+\frac{2b^2_{p,m}}{2p+5},&p'=p,\\
\frac{2b_{p,m}}{2p+5},&p'=p+2,\\
0&\mbox{otherwise},
\end{cases}
\end{split}
\end{equation}
Following along the similar lines as discussed in CQ and NP methods, the variational 
formulation works out to be
\begin{multline}
-\alpha_1^{-2}\left(\partial^2_{y_1}\wtilde{w}_m^{j+1},\wtilde{\theta}_{pm}\right)
+(1+\alpha_2^{-2}m^2)\left(\wtilde{w}_m^{j+1},\wtilde{\theta}_{pm}\right)\\
=(\wtilde{u}_m^j,\wtilde{\theta}_{pm}) +\mathcal{I}
\end{multline}
where
\begin{multline*}
\mathcal{I}= \alpha_1(1+\alpha_2^{-2}m^2)\left[ 
(\wtilde{\mathcal{B}}^{j+1}_{m,r}\chi_{m,r},\wtilde{\theta}_{pm})
-(\wtilde{\mathcal{B}}^{j+1}_{m,l}\chi_{m,l},\wtilde{\theta}_{pm})\right].
\end{multline*}
To further simplify the terms on right hand side, let $\wtilde{f}_m^j$ be such that
\begin{multline}
(1+\alpha_2^{-2}m^2)
[\alpha_1\wtilde{\mathcal{B}}^{j+1}_{m,r}\wtilde{\chi}_{m,r}(y_1)
-\alpha_1\wtilde{\mathcal{B}}^{j+1}_{m,l}\wtilde{\chi}_{m,l}(y_1)]\\
+\wtilde{u}^j_m=\wtilde{f}^j_m,
\end{multline}
then it follows that
\begin{multline}
\wtilde{f}^j_m =\wtilde{u}^j_m +\frac{\alpha_1\mathcal{D}_m}{2\kappa_{m}}
\left[\wtilde{\mathcal{B}}^{j+1}_{m,r}+\wtilde{\mathcal{B}}^{j+1}_{m,l}\right]L_0(y_1)\\
+\frac{\alpha_1\mathcal{D}_m}{2(1+\kappa_m)}
\left[\wtilde{\mathcal{B}}^{j+1}_{m,r}-\wtilde{\mathcal{B}}^{j+1}_{m,l}\right]L_1(y_1),
\end{multline}
where $ \mathcal{D}_m = (1+\alpha_2^{-2}m^2)$. 
Let $\wtilde{\wtilde{u}}_{p,m}$ denote the expansion coefficients of the field 
$\wtilde{{u}}_{m}$ in Legendre basis for each $m$-th Fourier component
and $\what{\wtilde{w}}_{p,m}$ denote the expansion coefficients in the boundary-adapted basis as
\begin{equation}
\begin{split}
&\wtilde{w}^{j+1}_m
= \sum_{p\in\field{J}_1}\what{\wtilde{w}}^{j+1}_{p,m}\theta_{p,m}(y_1)
= \sum_{p'\in\field{J}_1}\what{\wtilde{w}}^{j+1}_{p,m}\phi_{p,m}(y_1),\\
&\wtilde{u}_m^{j}
=\sum_{p=0}^N \wtilde{\wtilde{u}}^{j+1}_{p,m}L_{p}(y_1).
\end{split}
\end{equation}
The inner product $(\wtilde{f}_m^j,\wtilde{\theta}_{p,m})$ can be computed from the
knowledge of the Legendre coefficients of the field $\wtilde{f}_m^j$ using the
specific form of the boundary-adapted basis. This can be achieved via
\emph{quadrature matrix} given in~\ref{app:ip}.
The linear system thus becomes
\begin{equation*}
\left[{\alpha_1^{-2}}\wtilde{S}_{m,1}
+\left(1+{\alpha_2^{-2}}{m^2}\right)\wtilde{M}_{m,1}\right]
\what{\wtilde{\vv{w}}}^{\;j+1}_m
=\wtilde{Q}_{m,1}\Gamma\wtilde{\wtilde{\vv{f}}}^{\;j}_m.
\end{equation*}
Introducing the matrices
\begin{equation}
\left\{\begin{aligned}
&S_1=\diag\left(\left\{\wtilde{S}_{m,1}\right\}_{m\in\field{J}_2}\right),\\
&M_1=\diag\left(\left\{(1+\alpha_2^{-2}m^2)\wtilde{M}_{m,1}\right\}_{m\in\field{J}_2}\right),\\
&Q_1=\diag\left(\left\{\wtilde{Q}_{m,1}\Gamma\right\}_{m\in\field{J}_2}\right),
\end{aligned}\right.
\end{equation}
we can rewrite the linear system as
\begin{equation}
\left[\alpha_1^{-2}{S}_{1}+{M}_{1}\right]
\what{{\vv{w}}}^{j+1}
={Q}_{1}\wtilde{{\vv{f}}}^{\;j},
\end{equation}
where $\what{\vv{w}}$ and $\wtilde{\vv{f}}$ denote the block column vectors obtained 
by stacking the column vectors $\what{\wtilde{\vv{w}}}_m$ and $\wtilde{\wtilde{\vv{w}}}_m$, 
respectively, below one another for each $m$-th Fourier component. 
The linear system is then solved using LU-decomposition method. 
\section{Numerical Experiments--2D}\label{sec:numerical-experiments-2d}
In this section, we will carry out the extensive numerical tests to showcase the 
accuracy of the numerical schemes developed in this work to solve the 
IBVP~\eqref{eq:2D-SE-CT}, focusing on the case $\beta=+1$.
We start by analyzing the behaviour of the exact solutions
for the IBVP under consideration followed by studying the error evolution behaviour
and convergence analysis of the numerical schemes.

\subsection{Exact solutions}\label{sec:exact-solution}
The exact solutions admissible for our numerical experiments are such that the 
initial profile must be effectively supported within the computational domain.
We primarily consider the wavepackets which are a modulation of a Gaussian
envelop so that the requirement of effective initial support can be easily met.
In this class of solutions, we consider the Fourier-Chirped-Gaussian and the
Fourier-Hermite-Gaussian profiles.
\subsubsection{Fourier-Chirped-Gaussian profile}
\def\arraystretch{2}
\setlength{\tabcolsep}{1mm}
\begin{table*}[htb]
\centering
\caption{\label{tab:fcg2d} Fourier-Chirped-Gaussian profile with $A_0=2$ and 
$c_0\in\{4,8,12,16\}$.}
\begin{tabular}{m{15mm}m{10mm}m{20mm}m{20mm}m{25mm}m{25mm}}
\hline
\multicolumn{6}{c}{
$G\left(\vv{x},t\right) 
= A_0\sum_{j=1}^n G(x_1,t;{a}_j,{b}_j,{c}_j,\zeta_j),
\quad{c}_j=sgn(c_j)c_0 \quad \zeta_j=(\pi/d)K_j$ } \\
\hline
Type &$n$ & $({a}_j\in\field{R})_{j=1}^n$ & $({b}_j\in\field{R})_{j=1}^n$ 
& $(s_j=\sgn(c_j))_{j=1}^n$ & $(K_j\in\field{Z})_{j=1}^n$\\
\hline
I & $2$ & 
${a}_1=1/2.5$,\newline ${a}_2=1/2.3$
& $1/2$ &
${s}_1 = +1$,\newline 
${s}_2 = -1$
&
$K_1 = +2$,\newline 
$K_2 = -2$
\\
\hline
II& $4$ & 
${a}_1 = 1/2.5$,\newline 
${a}_2 = 1/2.3$,\newline
${a}_3 = 1/2.2$,\newline
${a}_4 = 1/2.4$
& $1/2$ &
${s}_1 = +1$,\newline 
${s}_2 = -1$,\newline
${s}_3 = +1$,\newline
${s}_4 = -1$ 
&
$K_1 = +2$,\newline 
$K_2 = -2$,\newline 
$K_3 = +4$,\newline 
$K_4 = -4$
\\
\hline
\end{tabular}
\end{table*}
Introduce the parameters $a\in\field{R}_+$ and $b\in\field{R}$ to 
define the function $\mathcal{G}=\mathcal{G}(x,t;a,b)$ by
\begin{equation}
\mathcal{G}=\frac{1}{\sqrt{1+4i(a+ib)t}}
\exp\left[\frac{-(a+ib)}{1+4i(a+ib)t}x^2\right].
\end{equation}
The family of solutions referred to as Fourier-chirped-Gaussian profile is then 
given by
\begin{multline}
G(\vv{x},t;a,b,c,\zeta)
=\mathcal{G}(x_1-c_1t,t;a,b)\\
\exp\left(+i\frac{1}{2}{c} x_1-i\frac{1}{4}{c}^2\,t\right)
\exp\left(i\zeta x_2 -i\zeta^2t\right),
\end{multline}
where ${a}\in\field{R}_+$ determines the effective support of the 
profile at $t=0$, ${b}\in\field{R}$ is the \emph{chirp} parameter, 
$\zeta \in\field{R}$ is the \emph{Fourier spectral} parameter and 
${c}\in\field{R}$ is the \emph{speed} of the profile. Using linear 
combination, one can further define a more general family of solutions with 
parameters $A_0,c_0\in\field{R}$,  
$({a}_j\in\field{R}_+)_{j=1}^n$, $({b}_j\in\field{R})_{j=1}^n$ and 
$(\zeta_j\in\field{R})_{j=1}^n$ given by
\begin{multline}
G\left(\vv{x},t;\;c_0,A_0;\;({a}_j)_{j=1}^n,({b}_j)_{j=1}^n,(\zeta_j)_{j=1}^n\right)\\ 
= A_0\sum_{j=1}^n G(\vv{x},t;{a}_j,{b}_j,{c}_j,\zeta_j),
\end{multline}
The specific values of the parameters of the solutions used in the
numerical experiments can be summarized as in Table~\ref{tab:fcg2d}. The energy 
content of the profile within the computational domain $\Omega_i$ 
over time is
\begin{equation}\label{eq:cg2d-energy-content}
E_{\Omega_i}(t)
=\left.{\int_{\Omega_i}|G(\vv{x},t)|^2d^2\vv{x}}\right/{\int_{\Omega_i}|G(\vv{x},0)|^2d^2\vv{x}}.
\end{equation}
The profiles are chosen with non-zero speed $c_0$ so that the field hits the
boundary of $\Omega_i$. 


\subsubsection{Fourier-Hermite-Gaussian profile}
\def\arraystretch{2}
\setlength{\tabcolsep}{1mm}
\begin{table*}[htb]
\centering
\caption{\label{tab:fhg2d} Fourier-Hermite-Gaussian profile with $A_0=2$ and 
$c_0\in\{4,8,12,16\}$.}
\begin{tabular}{m{15mm}m{10mm}m{20mm}m{20mm}m{25mm}m{25mm}}
\hline
\multicolumn{6}{c}{
$G\left(\vv{x},t\right) 
= A_0\sum_{j=1}^n G(x_1,t;{a}_j,{m}_j,{c}_j,\zeta_j),
\quad{c}_j=sgn(c_j)c_0 \quad \zeta_j=(\pi/d)K_j$ } \\
\hline
Type &$n$ & $({a}_j\in\field{R})_{j=1}^n$ & $({m}_j\in\field{N}_0)_{j=1}^n$ 
& $(s_j=\sgn(c_j))_{j=1}^n$ & $(K_j\in\field{Z})_{j=1}^n$\\
\hline
I & $2$ & 
${a}_1=1/2.5$,\newline ${a}_2=1/2.3$
& ${m}_1=1$,\newline ${m}_2=2$
& ${s}_1 = +1$,\newline 
${s}_2 = -1$
&
$K_1 = +2$,\newline 
$K_2 = -2$
\\
\hline
II& $4$ & 
${a}_1 = 1/2.5$,\newline 
${a}_2 = 1/2.3$,\newline
${a}_3 = 1/2.2$,\newline
${a}_4 = 1/2.4$
& 
${m}_1 = 1$,\newline 
${m}_2 = 2$,\newline
${m}_3 = 1$,\newline
${m}_4 = 2$
&
${s}_1 = +1$,\newline 
${s}_2 = -1$,\newline
${s}_3 = +1$,\newline
${s}_4 = -1$ 
&
$K_1 = +2$,\newline 
$K_2 = -2$,\newline 
$K_3 = +4$,\newline 
$K_4 = -4$
\\
\hline
\end{tabular}
\end{table*}

Consider the class of normalized Hermite-Gaussian functions defined by 
\begin{multline*}
\mathcal{G}_{m}(x,t;a)\\=\gamma_m^{-1}H_{m}\left(\frac{\sqrt{2a}\,x }{w(t)}\right)
\sqrt{\frac{\mu(t)}{a}} \exp\left[-\mu(t)x^2-i\,m\theta(t)\right],
\end{multline*}
where $a>0,\,m\in\field{N}_0$, $w(t)=\sqrt{1+(4at)^2}$ and 
\begin{equation*}
\frac{1}{\mu(t)}=\frac{1}{a}+i{4t}=\frac{1}{a}w(t)\exp[i\theta(t)],
\end{equation*}
with the normalization factor: 
$\gamma^2_{m}=2^m(m!)\sqrt{\pi}(2a)^{-1/2}$.

The Hermite polynomials are evaluated using the following relations
\begin{equation}
\begin{split}
& H_{n+1}(x)=2xH_n(x)-2nH_{n-1}(x),\\
&H_{n+1}'(x)=2(n+1)H_{n}(x),
\end{split}
\end{equation}
with $H_0(x) = 1$ and $H_1(x)=2x$.
Using these functions, we
can define a family of solutions referred to as Fourier-Hermite-Gaussian profile by
\begin{multline}
G(\vv{x},t;{m},{a},{c})=\mathcal{G}_{m_1}(x_1-c_1t,t;a_1)\\
\exp\left(+i\frac{1}{2}{c}{x}-i\frac{1}{4}c^2t\right)
\exp\left(i\zeta x_2 -i\zeta^2t\right),
\end{multline}
where ${m}\in\field{N}_0$ is the order parameter, 
${a}\in\field{R}_+$ determines the effective support of the profile 
at $t=0$, $\zeta \in\field{C}$ is the \emph{spectral} parameter and ${c}\in\field{R}$ 
is the velocity of the profile. Once again, using linear combination, one can further 
define a more general family of solutions with parameters $A_0,c_0\in\field{R}$,  
$({m}_j\in\field{N}_0)_{j=1}^n$, $(\zeta_j\in\field{C})_{j=1}^n$, and
$({a}_j\in\field{R}_+)_{j=1}^n$ given by
\begin{multline}
G\left(\vv{x},t;\;c_0,A_0;\;({m}_j)_{j=1}^n,({a}_j)_{j=1}^n,(\zeta_j)_{j=1}^n\right) \\
= A_0\sum_{j=1}^n G(\vv{x},t;{m}_j,{a}_j,{c}_j,\zeta_j).
\end{multline}
The specific values of the parameters of the solutions used in the
numerical experiments can be summarized as in Table~\ref{tab:fhg2d}. 
The profiles are chosen with non-zero speed $c_0$ so that the field hits the
boundary of $\Omega_i$. 
\begin{figure*}[!htb]
\begin{center}
\includegraphics[width=\textwidth]{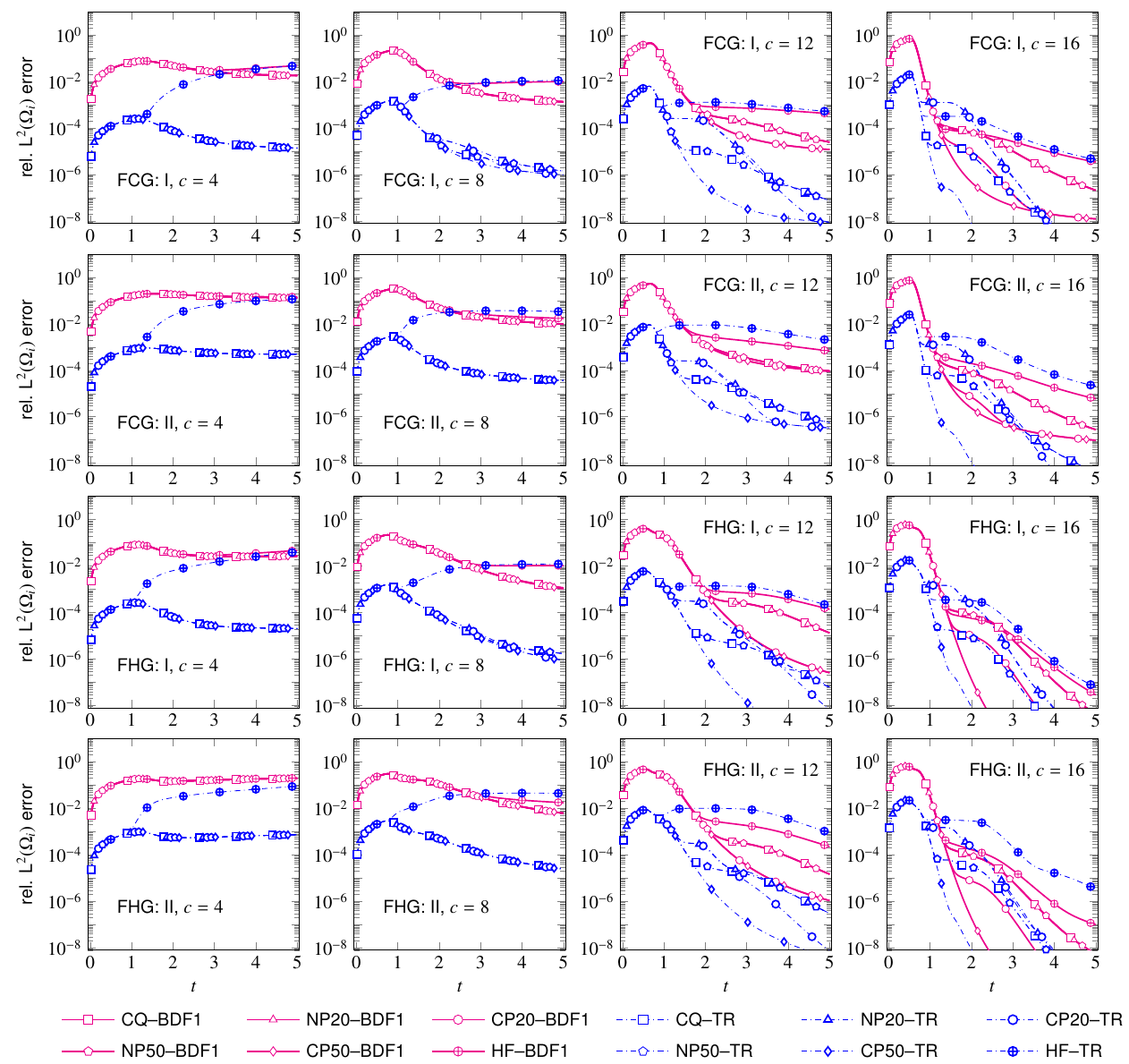}
\end{center}
\caption{\label{fig:wpfcgfhg-ee}The figure shows a comparison of evolution of error 
in the numerical solution of the IBVP~\eqref{eq:2D-SE-CT} with various approximations
of the TBCs corresponding to the Fourier-chirped-Gaussian and the Fourier-Hermite-Gaussian 
profiles with different values of the speed `c' (see Table~\ref{tab:fcg2d} and 
Table~\ref{tab:fhg2d}). The numerical parameters and the labels are described in 
Sec.~\ref{sec:tests-ee} where the error is quantified by~\eqref{eq:error-ibvp}.}
\end{figure*}
\def\arraystretch{2}
\setlength{\tabcolsep}{1mm}
\begin{table*}[htb]
\centering
\caption{\label{tab:ee-params}Numerical parameters for studying the evolution 
error}
\begin{tabular}{m{80mm}m{50mm}}\hline
Computational domain ($\Omega_i$) & $(-10,10)\times[-\pi,\pi)$\\\hline
Maximum time ($T_{max}$)          & $5$\\\hline
No. of time-steps ($N_t$)         & $5000+1$\\\hline
Time-step ($\Delta t$)            & $10^{-3}=T_{max}/(N_t-1)$\\\hline
\>Number of LGL-points ($(N+1)\times (N+1)$) & $200\times 200$\\\hline
\end{tabular}
\end{table*}
\subsection{Tests for evolution error}\label{sec:tests-ee}
In this section, we consider the IBVP in~\eqref{eq:2D-SE-CT} where the initial 
condition corresponds to the exact solutions described in Sec.~\ref{sec:exact-solution}. 
The error in the evolution of the profile computed numerically is quantified by the 
relative $\fs{L}^2(\Omega_i)$-norm: 
\begin{equation}\label{eq:error-ibvp}
e(t_j)=\frac{
\left(\int_{\Omega_i}
\left|u(\vv{x},t_j)-[u(\vv{x},t_j)]_{\text{num.}}
\right|^2d^2\vv{x}\right)^{1/2}\,}
{\left(\int_{\Omega_i}\left|u_0(\vv{x})\right|^2d^2\vv{x}\right)^{1/2}},
\end{equation}
for $t_j\in[0,T_{max}]$. The integral above will be approximated by Gauss 
quadrature over LGL-points. 

In this work, we have developed four discrete versions of the DtN maps, namely, 
CQ, NP, CP and HF. Each of these methods have a variant determined by the choice of 
one-step method used in the temporal discretization so that the complete list of
methods to be tested can be labelled as CQ--BDF1, NP--BDF1, CP--BDF1, HF--BDF1 
(corresponding to the one-step method BDF1), and, CQ--TR, NP--TR, CP--TR, HF--TR 
(corresponding to the one-step method TR). For the Pad\'e approximant based methods
`NP' and `CP', we distinguish the diagonal approximant of order $20$ from that of 
$50$ via the labels `NP20', `CP20'  and `NP50', `CP50'. The numerical parameters 
used in this section are summarized in Table~\ref{tab:ee-params}. The exact solutions 
used are described in Sec.~\ref{sec:exact-solution}.

The numerical results for the evolution error on $\Omega_i$ corresponding to the 
Fourier-chirped-Gaussian and the Fourier-Hermite-Gaussian profiles are shown in 
Fig.~\ref{fig:wpfcgfhg-ee}. It can be seen that the diagonal Pad\'e 
approximants based method, namely, NP and CP, perform equally well as compared 
to that of the convolution quadrature based method. 
Observe that the TR methods perform better than the BDF1 methods which is clear from the 
error peaks in Fig.~\ref{fig:wpfcgfhg-ee}. Further, note that the error curves for 
HF--TR starts separating from the rest just after certain initial time-steps. Also,
the performance of HF methods seems to improve for the case of faster moving profiles. 
For rest of the methods, the accuracy deteriorates with increase in the value of $c$
which can be seen from the error peaks in the Fig~~\ref{fig:wpfcgfhg-ee}. 
For faster moving profiles, the performance of CP50--TR turns out to be superior than the 
NP50--TR and CQ-TR. Also, NP20--TR and CP20--TR seems to be least accurate among TR 
methods which can be attributed to lower order of the Pad\'e approximant chosen. 
There is no significant difference in the results for evolution error on $\Omega_i$ 
corresponding to the Fourier-chirped-Gaussian profile and the Fourier-Hermite-Gaussian profile.

\def\arraystretch{2}
\setlength{\tabcolsep}{1mm}
\begin{table*}[htb]
\centering
\caption{\label{tab:ca-params}Numerical parameters for studying the convergence 
error}
\begin{tabular}{m{80mm}m{50mm}}\hline
Computational domain ($\Omega_i$) & $(-10,10)\times[-\pi,\pi)$\\\hline
Maximum time ($T_{max}$)          & $5$\\\hline
Set of no. of time-steps ($\field{N}_t$)& $\{2^8,2^9,\ldots,2^{16}\}$\\\hline
Time-step               & $\{\Delta t=T_{max}/(N_t-1),\;N_t\in\field{N}_t\}$\\\hline
\>Number of LGL-points ($(N+1)\times (N+1)$) & $200\times 200$\\\hline
\end{tabular}
\end{table*}
\subsection{Tests for convergence}\label{sec:tests-ca}
In this section, we analyze the convergence behaviour of the numerical schemes
for the IBVP in~\eqref{eq:2D-SE-CT} where the initial condition corresponds to the 
exact solutions described in Sec.~\ref{sec:exact-solution}. 
The error used to study the convergence behaviour is quantified by 
the maximum relative $\fs{L}^2(\Omega_i)$-norm: 
\begin{equation}\label{eq:max-error-ibvp}
e=\max\left\{e(t_j)|\;j=0,1,\ldots,N_t-1\right\},
\end{equation}
for $t_j\in[0,T_{max}]$. The Table~\ref{tab:ca-params} lists the numerical parameters 
for studying the convergence behaviour of CQ, NP, CP and HF methods. The numerical 
solution is labelled as described in Sec.~\ref{sec:tests-ee}. 

\begin{figure*}[!htb]
\begin{center}
\includegraphics[width=\textwidth]{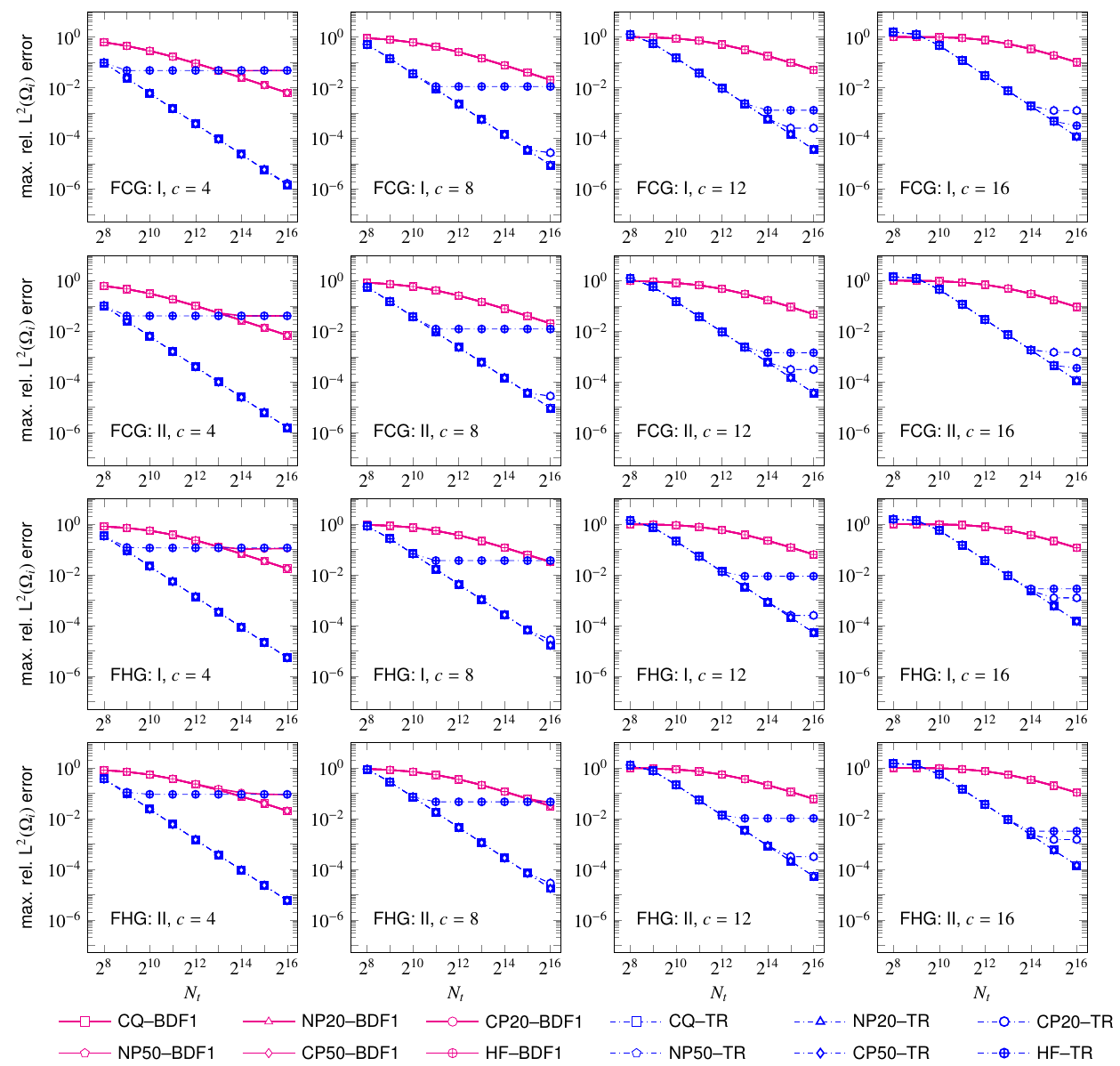}
\end{center}
\caption{\label{fig:wpfcgfhg-ca}The figure depicts the convergence behaviour of the
numerical solution of the IBVP~\eqref{eq:2D-SE-CT} with various approximations
of the TBCs corresponding to the Fourier-chirped-Gaussian and the Fourier-Hermite-Gaussian 
profiles with different values of the speed `c' (see Table~\ref{tab:fcg2d} and 
Table~\ref{tab:fhg2d}). The numerical parameters and the labels are described in 
Sec.~\ref{sec:tests-ca} where the error is quantified by~\eqref{eq:max-error-ibvp}.}
\end{figure*}

The numerical results for the convergence behaviour on $\Omega_i$ corresponding to the 
Fourier-chirped-Gaussian and the Fourier-Hermite-Gaussian profiles are shown in 
Fig.~\ref{fig:wpfcgfhg-ca}. It can be seen that all the methods under consideration
for the convergence tests show stable behaviour for the parameters specified in 
Table~\ref{tab:ca-params}. In the log-log scale, we can clearly identify the error 
curve to be a straight line before it plateaus. The empirical orders can be recovered 
from the slope which we found consistent with the order of the underlying one-step method. 
The TR methods perform better than the BDF1 methods which is obvious from the 
slope of the error curves in Fig.~\ref{fig:wpfcgfhg-ca}. Particularly for the HF methods,
we note that the error background decreases with increase in the value of $c$ as one 
would expect with the case of high-frequency approximation of the TBCs.  
Also observe that CP20 methods starts plateauing after a certain number of time-steps 
unlike CP50. This can be attributed to the approximate nature of the TBCs with a lower 
order of diagonal approximants chosen. 
Finally, let us observe that there is no significant difference in the 
results for convergence on $\Omega_i$ corresponding to the 
Fourier-chirped-Gaussian profile and the Fourier-Hermite-Gaussian profile.  
\section{Numerical Implementation--3D}\label{sec:numerical-implementation-3d}
In this section, we address the complete numerical solution 
of the initial boundary-value problem (IBVP) stated in~\eqref{eq:3D-SE-CT}.
The formulation of the discrete linear system for the IBVP requires 
temporal as well as spatial discretization of the problem. 
For the computational domain $\Omega_i$, we introduce a reference domain 
$\Omega_i^{\text{ref.}}=\field{I}\times\pi\field{I}\times\pi\field{I}$ keeping in 
view the typical domain of definition of the orthogonal polynomials being used.
In order to describe the associated linear maps between the reference domain and 
the actual computational domain, we introduce the variables 
$y_1,y_2,y_3\in\Omega_i^{\text{ref.}}$ such that
\begin{equation*}
\begin{split}
&x_1 = J_1 y_1+\bar{x}_1,\quad J_1 = \frac{1}{2}(x_r-x_l),\quad
\bar{x}_1=\frac{1}{2}(x_l+x_r), \\
&x_2 = J_2 y_2,\quad J_2 = \frac{d_2}{\pi},\quad \beta_1=J_1^{-2},\quad
\beta_2=\beta J_2^{-2},\\
&x_3 = J_3 y_3,\quad J_3 = \frac{d_3}{\pi},\quad \beta_3=\beta J_3^{-2},
\end{split}
\end{equation*}
Let $\Delta t$ denote the time-step. 
In the rest of the section, with a slight abuse of notation, we switch to the 
variables $y_1\in\field{I}$ and $y_2=y_3\in\pi\field{I}$ for all the discrete approximations 
of the dependent variables. For instance, $u^j(\vv{y})$ is taken to approximate 
$u(\vv{x},j\Delta t)$ for $j=0,1,2,\ldots,N_t-1$. The temporal discretization of the 
interior problem using one-step methods is enumerated below:
\begin{itemize}
\item The BDF1 based discretization is given by
\begin{equation}\label{eq:3d-se-bdf1}
\begin{split}
&i\frac{u^{j+1}-u^{j}}{\Delta t}+\triangle u^{j+1}=0, \quad \rho = 1/\Delta t,\\
&\left(\beta_1\partial^2_{y_1}
+\beta_2\partial^2_{y_2}+\beta_3\partial^2_{y_3}\right)u^{j+1}+i\rho u^{j+1}
=i\rho u^j.
\end{split}
\end{equation}
\item The TR based discretization is given by
\begin{equation}\label{eq:3d-disc-tr}
\begin{split}
&i\frac{u^{j+1}-u^{j}}{\Delta t}+\triangle u^{j+1/2}=0,\quad \rho = 2/\Delta t,\\
&\left(\beta_1\partial^2_{y_1}+\beta_2\partial^2_{y_2}
+\beta_3\partial^2_{y_3}\right)v^{j+1}+i\rho v^{j+1}
=i\rho u^j,
\end{split}
\end{equation}
\end{itemize}
where we have used the staggered samples of the field in the last step.
The rest of this section is organized as follows: Sec.~\ref{sec:dtn-maps-3d}
addresses the temporal discretization of the boundary maps which is followed 
by a discussion of the resulting spatial problem in Sec.~\ref{sec:IBVP-3d}.
\subsection{Discretizing the boundary conditions}\label{sec:dtn-maps-3d}
As discussed above, we aim at a compatible temporal discretization of the novel 
boundary conditions formulated in Sec.~\ref{sec:CT-NP-3D} in order to solve the 
IBVP in~\eqref{eq:3D-SE-CT}.
In order to facilitate the temporal discretization of the TBCs, we introduce the 
equispaced samples of the auxiliary function, $\varphi^{j,k}(\vv{y})$, to 
approximate $\varphi(\vv{x},j\Delta t,k\Delta t)$ for $j,k=0,1,\ldots N_t-1$. It 
follows from the definition of the auxiliary function that the diagonal values of 
the auxiliary function at the discrete level are determined by the interior field,
i.e., $\varphi^{j,j}(y_1,y_2,y_3) = u^j(y_1,y_2,y_3)$.
For the space discretization of the auxiliary equation on the faces of $\Omega_i$, 
we use the Fourier-Galerkin method which employs the series
\begin{equation}
\varphi^{j,k}(\vv{y})
=\sum_{m_2\in\field{J}_2}\sum_{m_3\in\field{J}_3}\wtilde{\varphi}_{m_2,m_3}^{j,k}(\vv{y})
e^{i(m_2y_2+m_3y_3)},
\end{equation}
where the index set is given by
\begin{equation}
\field{J}_p=\left\{-\frac{N_p}{2},\ldots,\frac{N_p}{2}-1\right\},\quad p=2,3.
\end{equation}
For the sake of brevity of presentation, let
$\wtilde{\varphi}_{\vv{m}}^{j,k}$ represent the Fourier samples of the auxiliary
function i.e. $\wtilde{\varphi}_{m_2,m_3}^{j,k}$.
\subsubsection{NP-BDF1}
The discrete scheme for the boundary maps is said to be `NP--BDF1' if the 
temporal derivatives in the realization of the boundary maps, dubbed as `NP', 
developed in Sec.~\ref{sec:CT-NP-3D} are discretized using the one-step method 
BDF1. To achieve this, we start by writing the time-discrete form of the DtN 
maps present in~\eqref{eq:maps-npade-3d} in the reference domain as
\begin{equation}\label{eq:map-npade-bdf1-3d}
\sqrt{\beta_1}\partial_{y_1}\wtilde{u}_{\vv{m}}^{j+1}
\pm e^{-i\pi/4}\left[ b_0\wtilde{u}_{\vv{m}}^{j+1} -
\sum_{k=1}^M b_k\wtilde{\varphi}^{j+1,j+1}_{k,\vv{m}}\right]=0.
\end{equation}
In order to turn these equations into a Robin-type boundary condition, we need to 
compute the discrete samples $\varphi^{j+1,j+1}_{k,\vv{m}}$ 
by considering the ODEs established for these auxiliary function earlier. 
Here, we discuss the BDF1-based discretization of~\eqref{eq:auxi-npade-3d} by first
employing Fourier expansion: 
Setting $\rho=1/\Delta t$ and for $k=1,2,\ldots,M$, we have
\begin{multline}
\frac{\wtilde{\varphi}^{j+1,j+1}_{k,\vv{m}}-\wtilde{\varphi}^{j, j+1}_{k,\vv{m}}}{\Delta t} 
+\eta^2_k\wtilde{\varphi}^{j+1,j+1}_{k,\vv{m}}
=\wtilde{\varphi}_{\vv{m}}^{j+1,j+1}\\
\wtilde{\varphi}^{j+1,j+1}_{k,\vv{m}} 
=\frac{\rho}{(\rho+\eta^2_k)}\wtilde{\varphi}^{j,j+1}_{k,\vv{m}}
+\frac{1}{(\rho+\eta^2_k)}\wtilde{u}_{\vv{m}}^{j+1}.
\end{multline}
Introduce the scaling $\ovl{\eta}_k=\eta_k/\sqrt{\rho}$ so that
\begin{equation}
\wtilde{\varphi}^{j+1,j+1}_{k,\vv{m}} 
=\frac{1}{(1+\ovl{\eta}^2_k)}\left[\wtilde{\varphi}^{j,j+1}_{k,\vv{m}}
+\frac{1}{\rho}\wtilde{u}_{\vv{m}}^{j+1}\right].
\end{equation}
The value of $\wtilde{\varphi}^{j,j+1}_{k,\vv{m}}$ can be obtained by propagation 
along $\tau_{\perp}$ using $\wtilde{\varphi}^{j,j}_{k,\vv{m}}$ as the initial condition 
in~\eqref{eq:FSE-auxi-ivp-3d}:
\begin{equation}
\wtilde{\varphi}^{j,j+1}_{k,\vv{m}}
=\frac{1}{(1+\alpha_2^{-2}m_2^2+\alpha_3^{-2}m_3^2)}\wtilde{\varphi}^{j,j}_{k,\vv{m}},
\end{equation}
where $m_p\in\field{J}_p,\;p=2,3$.
Next we plug-in the value of the auxiliary field,
$\wtilde{\varphi}^{j+1,j+1}_{k,\vv{m}}$,
in~\eqref{eq:map-npade-bdf1-3d}
to obtain
\begin{multline}
\partial_{y_1}\wtilde{u}_{\vv{m}}^{j+1}\pm\alpha_1\varpi\wtilde{u}_{\vv{m}}^{j+1}\\
=\mp\frac{\alpha_1}{(1+\alpha_2^{-2}m_2^2+\alpha_3^{-2} m_3^2)}\left[\sum_{k=1}^M 
\Gamma_k\wtilde{\varphi}^{j,j}_{k,\vv{m}}\right]\\
=\mp\alpha_1\wtilde{\mathcal{B}}^{j+1}_{\vv{m},a_1},
\end{multline}
where we have set $\ovl{b}_0=b_0/\sqrt{\rho},\; \ovl{b}_k =b_k/\sqrt{\rho}$ and
\begin{equation}
\varpi =\ovl{b}_0 +\frac{1}{\rho}\sum_{k=1}^M\Gamma_k,\quad
\Gamma_{k} =-{\ovl{b}_{k}}/{(1+\ovl{\eta}^2_{k})}.
\end{equation}
Note that these Robin type maps are written for each component in the Fourier
expansion. 
\subsubsection{NP--TR}
The discrete scheme for the boundary maps is said to be `NP--TR' if the 
temporal derivatives in the realization of the boundary maps, dubbed as `NP', 
developed in Sec.~\ref{sec:CT-NP-3D} are discretized using the one-step method 
TR. To achieve this, we start by writing the time-discrete form of the DtN 
maps present in~\eqref{eq:maps-npade-3d} in the reference domain as
\begin{equation}\label{eq:map-npade-tr-3d}
\sqrt{\beta_1}\partial_{y_1}\wtilde{v}_{\vv{m}}^{j+1}
\pm e^{-i\pi/4}\left[ b_0\wtilde{v}_{\vv{m}}^{j+1} -
\sum_{k=1}^M b_k\wtilde{\varphi}^{j+1/2,j+1/2}_{k,\vv{m}}\right]=0.
\end{equation}
In order to turn these equations into a Robin-type boundary condition, we need to 
compute the discrete samples $\varphi^{j+1/2,j+1/2}_{k,\vv{m}}$ 
by considering the ODEs established for these auxiliary function earlier. 
Here, we discuss the TR-based discretization of~\eqref{eq:auxi-npade-3d} by first
employing Fourier expansion: 
Setting $\rho=2/\Delta t$ and for $k=1,2,\ldots,M$, we have
\begin{multline}
\frac{\wtilde{\varphi}^{j+1,j+1}_{k,\vv{m}}-\wtilde{\varphi}^{j,j+1}_{k,\vv{m}}}{\Delta t} 
+\eta^2_k\frac{\wtilde{\varphi}^{j+1,j+1}_{k,\vv{m}}+\wtilde{\varphi}^{j,j+1,j}_{k,\vv{m}}}{2}\\
= \frac{\wtilde{\varphi}_{\vv{m}}^{j+1,j+1}+\wtilde{\varphi}_{\vv{m}}^{j,j+1}}{2}.
\end{multline}
Introduce the scaling $\ovl{\eta}_k=\eta_k/\sqrt{\rho}$ so that
\begin{multline}
\wtilde{\varphi}^{j+1,j+1}_{k,\vv{m}} 
=\frac{(1-\ovl{\eta}^2_k)}{(1+\ovl{\eta}^2_k)}\wtilde{\varphi}^{j,j+1}_{k,\vv{m}}\\
+\frac{2/\rho}{(1+\ovl{\eta}^2_k)}\left[\wtilde{v}_{\vv{m}}^{j+1}
+\frac{\wtilde{\varphi}_{\vv{m}}^{j,j+1}-\wtilde{\varphi}_{\vv{m}}^{j,j}}{2}\right].
\end{multline}
For the staggered configuration, we have
\begin{multline}\label{eq:map-novel-pade-trap-3d}
\wtilde{\varphi}^{j+1/2,j+1/2}_{k,\vv{m}} 
=\frac{1}{2}\left[\left(\frac{1-\ovl{\eta}^2_k}{1+\ovl{\eta}^2_k}\right)
\wtilde{\varphi}^{j,j+1}_{k,\vv{m}}+\wtilde{\varphi}^{j,j}_{k,\vv{m}}\right]\\
+\frac{1/\rho}{(1+\ovl{\eta}^2_k)}\left[\wtilde{v}_{\vv{m}}^{j+1}
+\frac{\wtilde{\varphi}_m^{j,j+1}-\wtilde{\varphi}_{\vv{m}}^{j,j}}{2}\right].
\end{multline}
The values of $\wtilde{\varphi}^{j,j+1}_{k,\vv{m}}$ and
$\wtilde{\varphi}^{j,j+1}_{\vv{m}}$ 
can be obtained by propagation along $\tau_{\perp}$ using
$\wtilde{\varphi}^{j,j}_{k,\vv{m}}$ and $\wtilde{\varphi}^{j,j}_{\vv{m}}$ as 
the initial condition:
\begin{multline}
\wtilde{\varphi}^{j,j+1}_{k,\vv{m}}
=\frac{\left(1-\alpha_2^{-2}m_2^2-\alpha_3^{-2}m_3^2\right)}{\left(1+\alpha_2^{-2}m_2^2
+\alpha_3^{-2}m_3^2\right)}
\wtilde{\varphi}^{j,j}_{k,\vv{m}},\\
\frac{1}{2}\left[\wtilde{\varphi}^{j,j+1}_{\vv{m}}-\wtilde{\varphi}^{j,j}_{\vv{m}}\right]
=\frac{\left(-\alpha_2^{-2}m_2^2-\alpha_3^{-2}m_3^2\right)}{\left(1+\alpha_2^{-2}m_2^2
+\alpha_3^{-2}m_3^2\right)}
\wtilde{\varphi}^{j,j}_{\vv{m}}.
\end{multline}
Next we plug-in the value of $\wtilde{\varphi}^{j+1/2,j+1/2}_{k,\vv{m}}$
in~\eqref{eq:map-npade-tr-3d} which leads to following compact form of the 
boundary maps
\begin{equation}
\partial_{y_1}\wtilde{v}_{\vv{m}}^{j+1}\pm\alpha_1\varpi\wtilde{v}_{\vv{m}}^{j+1}
=\mp\alpha_1\wtilde{\mathcal{B}}^{j+1/2}_{\vv{m},a_1},\\
\end{equation}
where $\varpi$ is same as that in BDF1 and the history functions are given by
\begin{equation}
\begin{split}
&\wtilde{\mathcal{B}}^{j+1/2}_{\vv{m},a_1}
=\sum_{k=1}^M\frac{-\ovl{b}_k}{2}\wtilde{\varphi}^{j,j}_{k,\vv{m}}\\
&+\sum_{k=1}^M\frac{\Gamma_k}{\rho}\left(
\frac{-\alpha_2^{-2}m_2^2-\alpha_3^{-2}m_3^2}{1+\alpha_2^{-2}m_2^2+\alpha_3^{-2}m_3^2}
\right)\wtilde{u}^{j}_{\vv{m}}\\
&+\sum_{k=1}^M\frac{-\ovl{b}_k}{2}\left(\frac{1-\ovl{\eta}^2_k}{1+\ovl{\eta}^2_k}\right)
\left(\frac{1-\alpha_2^{-2}m_2^2-\alpha_3^{-2}m_3^2}{1+\alpha_2^{-2}m_2^2+\alpha_3^{-2}m_3^2}\right)
\wtilde{\varphi}^{j,j}_{k,\vv{m}}\\
\end{split}
\end{equation}
Note that these Robin type maps are written for each component in Fourier
expansion. 
\def\arraystretch{2}
\setlength{\tabcolsep}{1mm}
\begin{table*}[htb]
\centering
\caption{\label{tab:ee-params-3d}Numerical parameters for studying the evolution 
error}
\begin{tabular}{m{80mm}m{50mm}}\hline
Computational domain ($\Omega_i$) & $(-10,10)\times[-\pi,\pi)\times[-\pi,\pi)$\\\hline
Maximum time ($T_{max}$)          & $5$\\\hline
No. of time-steps ($N_t$)         & $5000+1$\\\hline
Time-step ($\Delta t$)            & $10^{-3}=T_{max}/(N_t-1)$\\\hline
\>Number of LGL-points ($(N+1)\times (N+1)\times (N+1)$) & $100\times 100\times 100$\\\hline
\end{tabular}
\end{table*}
\subsection{Numerical Solution of the IBVP}\label{sec:IBVP-3d}
Having discussed the temporal discretization of the 3D problem, we now aim at
the spatial discretization using Legendre Galerin spectral method. Recall the form 
of the TBCs (for NP methods) casted as Robin-type boundary maps specific to 
each Fourier component:
\begin{equation}
\partial_{y_1}\wtilde{u}_{\vv{m}}^{j+1}\pm \alpha_1\varpi\wtilde{u}_{\vv{m}}^{j+1}
=\mp\alpha_1\wtilde{\mathcal{B}}^{j+1}_{\vv{m},a_1},\quad a_1\in\{r,l\}.
\end{equation}
Following along the similar lines as discussed for the case of 2D problem, we
define the space of boundary-adapted basis to enforce these maps exactly as
\begin{equation*}
\fs{W}_N=\left\{\wtilde{u}_{\vv{m}}\in\fs{P}_N\left|\;
\begin{aligned}
&(\partial_{y_1}-\kappa)\wtilde{u}_{\vv{m}}(y_1)|_{y_1=-1}=0,\\
&(\partial_{y_1}+\kappa)\wtilde{u}_{\vv{m}}(y_1)|_{y_1=+1}=0
\end{aligned}\right.
\right\},
\end{equation*}
where $\kappa=\alpha_1\varpi$ which is determined by the Robin-type formulation
of the TBCs for NP methods.
We introduce a function $\chi_{\vv{m}}$ to convert the discrete BCs to a homogeneous 
form by stating the original field as
$\wtilde{u}_{\vv{m}}^j=\wtilde{w}_{\vv{m}}^j+\chi_{\vv{m}}^{j}(y_1)$
where $\wtilde{w}_{\vv{m}}^j\in\fs{W}_N$ so that
\begin{equation*}
\begin{split}
[-\alpha_1^{-2}\partial^2_{y_1} +\alpha_2^{-2}m_2^2+\alpha_3^{-2}m_3^2 +1]
(\wtilde{w}_{\vv{m}}^{j+1}+\chi_{\vv{m}}^{j+1})=\wtilde{u}_{\vv{m}}^j,\\
\begin{aligned}
&\left.(\partial_{y_1}-\kappa)\wtilde{w}_{\vv{m}}^{j+1}\right|_{y_1=-1}=0,\\
&\left.(\partial_{y_1}+\kappa)\wtilde{w}_{\vv{m}}^{j+1}\right|_{y_1=+1}=0,
\end{aligned}
\end{split}
\end{equation*}
and the field $\chi_{\vv{m}}^j(y_1)$ is forced to satisfy the constraints:
\begin{equation}
\begin{split}
&(\partial_{y_1}-\kappa)\chi_{\vv{m}}^{j}(y_1)|_{y_1=-1}=+\alpha_1\wtilde{\mathcal{B}}^j_{{\vv{m}},l},\\
&(\partial_{y_1}+\kappa)\chi_{\vv{m}}^{j}(y_1)|_{y_1=+1}=-\alpha_1\wtilde{\mathcal{B}}^j_{{\vv{m}},r},
\end{split}
\end{equation}
Following in the similar manner as in 2D case, the lifting functions for
this case works out to be
\begin{equation}
\left\{\begin{aligned}
\chi_{l}(y_1) &=-\frac{1}{2\kappa}L_0(y_1)+\frac{1}{2(\kappa+1)}L_1(y_1),\\
\chi_{r}(y_1) &=+\frac{1}{2\kappa}L_0(y_1)+\frac{1}{2(\kappa+1)}L_1(y_1).
\end{aligned}\right.
\end{equation}
The lifted field becomes
\begin{equation}\label{eq:lifted-field-m-3d}
\wtilde{u}_{\vv{m}}^j({y_1})=\wtilde{w}_{\vv{m}}^j({y_1})
+\alpha_1\chi_{l}\wtilde{\mathcal{B}}^j_{{\vv{m}},l}-\alpha_1\chi_{r}\wtilde{\mathcal{B}}^j_{{\vv{m}},r},
\end{equation}
Next, we choose the basis $\theta_{p}\in\fs{W}_N$ defined as
\begin{equation}
\theta_{p}(y_1)=\phi_p(y_1),\quad p\in\field{J}_1,
\end{equation}
where $\{\phi_{p}|\,p\in\field{J}_1\}$ can be constructed using the ansatz 
$\phi_{p}(y_1)=L_p(y_1)+a_pL_{p+1}(y_1)+b_pL_{p+2}(y_1)$.
Imposing the boundary conditions in the definition of $\fs{W}_N$, the 
sequences $(a_p)_{p\in\field{J}_1}$ and $(b_p)_{p\in\field{J}_1}$ work out to be
\begin{equation}
a_p=0,\quad
b_p=-\frac{\kappa+\frac{1}{2}p(p+1)}{\kappa+\frac{1}{2}(p+2)(p+3)},
\quad p\in\field{J}_1.
\end{equation}
The mass and stiffness matrices with the boundary-adapted basis functions remain same as 
described for NP method (in 2D problem) due to the fact that $\varpi$ present
in TBCs is independent of the Fourier component.
The variational formulation then works out to be
\begin{multline}
-\alpha_1^{-2}\left(\partial^2_{y_1}\wtilde{w}_{\vv{m}}^{j+1},\theta_{p}\right)
+[1+\alpha_2^{-2}m_2^2+\alpha_3^{-2}m_3^2]\left(\wtilde{w}_{\vv{m}}^{j+1},\theta_{p}\right)\\
=(\wtilde{u}_{\vv{m}}^j,\theta_p) +\mathcal{I}
\end{multline}
where
\begin{equation}
\mathcal{I}= \alpha_1\mathcal{D}_{\vv{m}}\left[ 
(\wtilde{\mathcal{B}}^{j+1}_{{\vv{m}},r}\chi_{r},\theta_{p})
-(\wtilde{\mathcal{B}}^{j+1}_{{\vv{m}},l}\chi_{l},\theta_{p})\right].
\end{equation}
and $ \mathcal{D}_{\vv{m}} = (1+\alpha_2^{-2}m_2^2+\alpha_3^{-2}m_3^2)$.
To further simplify the terms on right hand side, let $\wtilde{f}_{\vv{m}}^j$ be such that
\begin{multline}
\mathcal{D}_{\vv{m}}
[\alpha_1\wtilde{\mathcal{B}}^{j+1}_{{\vv{m}},r}\wtilde{\chi}_{r}(y_1)
-\alpha_1\wtilde{\mathcal{B}}^{j+1}_{{\vv{m}},l}\wtilde{\chi}_{l}(y_1)]
+\wtilde{u}^j_{\vv{m}}=\wtilde{f}^j_{\vv{m}},
\end{multline}
then it follows that
\begin{multline}
\wtilde{f}^j_{\vv{m}} =\wtilde{u}^j_{\vv{m}} +\frac{\alpha_1\mathcal{D}_{\vv{m}}}{2\kappa}
\left[\wtilde{\mathcal{B}}^{j+1}_{{\vv{m}},r}+\wtilde{\mathcal{B}}^{j+1}_{{\vv{m}},l}\right]L_0(y_1)\\
+\frac{\alpha_1\mathcal{D}_{\vv{m}}}{2(1+\kappa)}
\left[\wtilde{\mathcal{B}}^{j+1}_{{\vv{m}},r}-\wtilde{\mathcal{B}}^{j+1}_{{\vv{m}},l}\right]L_1(y_1).
\end{multline}
Let $\wtilde{\wtilde{u}}_{p,{\vv{m}}}$ denote the expansion coefficients of the field 
$\wtilde{{u}}_{{\vv{m}}}$ in Legendre basis for each $\vv{m}$-th Fourier component
and $\what{\wtilde{w}}_{p,{\vv{m}}}$ denote the expansion coefficients in the boundary-adapted basis as
\begin{equation}
\begin{split}
&\wtilde{w}^{j+1}_{\vv{m}}
= \sum_{p\in\field{J}_1}\what{\wtilde{w}}^{\;j+1}_{p,{\vv{m}}}\theta_{p}(y_1)
= \sum_{p'\in\field{J}_1}\what{\wtilde{w}}^{\;j+1}_{p,{\vv{m}}}\phi_{p}(y_1),\\
&\wtilde{u}_{\vv{m}}^{j}
=\sum_{p=0}^N \wtilde{\wtilde{u}}^{\;j+1}_{p,{\vv{m}}}L_{p}(y_1).
\end{split}
\end{equation}
The inner product $(\wtilde{f}_{\vv{m}}^j,{\theta}_{p})$ can be computed from the
knowledge of the Legendre coefficients of the field $\wtilde{f}_{\vv{m}}^j$ using the
specific form of the boundary-adapted basis. This can be achieved via
\emph{quadrature matrix} given in~\ref{app:ip}.
The linear system thus becomes
\begin{equation*}
\left[{\alpha_1^{-2}}S_1
+\mathcal{D}_{\vv{m}}M_1\right]
\what{\wtilde{\vv{w}}}^{\;j+1}_{\vv{m}}
=Q_1\Gamma\wtilde{\wtilde{\vv{f}}}^{\;j}_{\vv{m}}.
\end{equation*}
The mass and stiffness matrices are defined in~\eqref{eq:sys-mat-2d} and the 
linear system is then solved using LU-decomposition method. 

\begin{figure*}[!htb]
\begin{center}
\includegraphics[width=\textwidth]{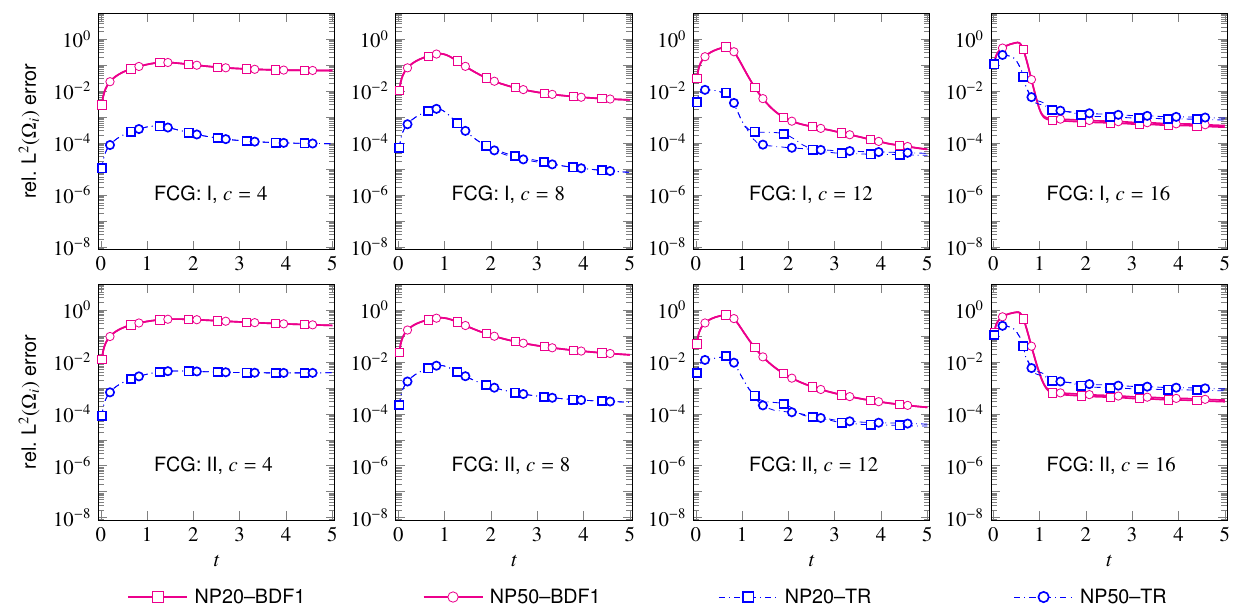}
\end{center}
\caption{\label{fig:wpf2cg-ee}The figure shows a comparison of evolution of error 
in the numerical solution of the IBVP~\eqref{eq:3D-SE-CT} with novel Pad\'e
based approximation of the TBCs corresponding to the Fourier-chirped-Gaussian 
profile with different values of the speed `c' (see Table~\ref{tab:fcg2d}). The 
numerical parameters and the labels are described in 
Sec.~\ref{sec:tests-ee-3d} where the error is quantified by~\eqref{eq:error-ibvp}.}
\end{figure*}
\section{Numerical Experiments--3D}\label{sec:numerical-experiments-3d}
In this section, we will carry out some numerical tests to showcase the 
accuracy of the numerical schemes developed in this work to solve the 
IBVP~\eqref{eq:3D-SE-CT}, focusing on the case $\beta=+1$.
We start by writing the exact solutions
for the IBVP under consideration followed by studying the error evolution behaviour
of the numerical schemes.
\subsection{Fourier-Chirped-Gaussian profile}\label{sec:exact-solution-3d}
For the 3D problem, we restrict ourselves to the Fourier-chirped-Gaussian
profiles for the choice of exact solutions. Recalling the function $\mathcal{G}$ 
from Sec.~\ref{sec:exact-solution}, we can define a family of solutions 
(for the 3D problem) referred to as Fourier-chirped-Gaussian 
profile given by
\begin{multline}
G(\vv{x},t;a,b,c,\vv{\zeta}_{\perp})
=\mathcal{G}(x_1-c_1t,t;a,b)\\
\exp\left(+i\frac{1}{2}{c} x_1-i\frac{1}{4}{c}^2\,t\right)
\exp\left(i\vv{\zeta}_{\perp}.\vv{x}_{\perp} -i\vv{\zeta}_{\perp}^2t\right),
\end{multline}
where $\vs{\zeta}_{\perp}$ denote the covariable corresponding to $\vv{x}_{\perp}$ 
used in the two dimensional Fourier transform. Rest of the parameters are same
as described in Sec.~\ref{sec:exact-solution}.
The specific values of the parameters of the solutions used in the
numerical experiments are summarized in Table~\ref{tab:fcg2d} assuming
$\zeta_2=\zeta_3=\zeta$ which means we take the periodicity to be same along $\vv{e}_2$
and $\vv{e}_3$.
The profiles are chosen with non-zero speed $c_0$ so that the field hits the
boundary of $\Omega_i$. 

\subsection{Tests for evolution error}\label{sec:tests-ee-3d}
In this section, we consider the IBVP in~\eqref{eq:3D-SE-CT} where the initial 
condition corresponds to the exact solutions described in
Sec.~\ref{sec:exact-solution-3d}. 
The error in the evolution of the profile computed numerically is quantified by the 
relative $\fs{L}^2(\Omega_i)$-norm defined in Sec.~\ref{sec:tests-ee}.

For the 3D problem, we only consider the novel-Pad\'e based discrete versions of 
the DtN maps, namely, NP. They have a variant determined by the choice of 
one-step method used in the temporal discretization so that the complete list of
methods to be tested can be labelled as NP--BDF1, 
(corresponding to the one-step method BDF1), and, NP--TR
(corresponding to the one-step method TR). For the Pad\'e approximant based 
method (NP), we distinguish the diagonal approximant of order $20$ from that of 
$50$ via the labels `NP20' and `NP50'. The numerical parameters 
used in this section are summarized in Table~\ref{tab:ee-params-3d}. The exact solutions 
used are defined in Sec.~\ref{sec:exact-solution-3d}. 

The numerical results for the evolution error on $\Omega_i$ corresponding to the 
Fourier-chirped-Gaussian profile are shown in Fig.~\ref{fig:wpf2cg-ee}. It can be 
seen that the diagonal Pad\'e approximants based method, namely, NP, show stable behaviour.  
Observe that the TR methods perform better than the BDF1 methods which is clear from the 
error peaks in Fig.~\ref{fig:wpf2cg-ee}. Note that the accuracy deteriorates with 
increase in the value of $c$ which can be seen from the error peaks in the 
Fig~~\ref{fig:wpf2cg-ee}. This behaviour can be attributed to the increased
oscillatory nature of the solution. Moreover, for the faster moving profiles,
error curves plateaus after certain number of time steps.

\section{Conclusion}\label{sec:conclusion}
In this work, we presented an efficient numerical implementation of the 
transparent boundary operator and its various approximations for the case of 
free Schr\"{o}dinger equation on a rectangular computational domain with
periodicity assumption along a subset of the unbounded directions. 
More specifically, we compared the performance of various approximations of the
transparent boundary operator, namely, convolution quadrature based approach, 
conventional Pad\'e approach, novel Pad\'e approach and high-frequency 
approximation method. Next, keeping in mind the superior accuracy and lower 
computational complexity of novel Pad\'e method (in comparison to rest of the 
approaches), we extend the novel Pad\'e approach for the TBCs
corresponding to free Schr\'odinger equation on a hyperrectangular computational
domain in 3D.

For the spatial problem, we use a Legendre-Fourier Galerkin method where a 
boundary adapted basis is constructed using Legendre polynomials that ensures 
the bandedness of the resulting linear system. Temporal discretization is 
addressed with two one-step methods, namely BDF1 and TR. As far as time-marching
is concerned, we would like to mention that any higher order one-step method can
be readily used with the novel Pad\'e method to further enhance the accuracy of the
numerical schemes. Several numerical tests
are carried out to demonstrate the efficiency and empirical stability of various 
numerical schemes for the 2D as well as the 3D problem under consideration.
Theoretical proofs for stability and convergence of our numerical schemes are
beyond the scope of this paper. We defer these issues to a future publication. 
Extension of the ideas presented in this paper to more general systems will be 
explored in future publications.

\section*{Acknowledgements}
The first author thanks CSIR India (grant no. 09/086(1431)/2019-EMR-I) for 
providing the financial assistance.  


\providecommand{\noopsort}[1]{}\providecommand{\singleletter}[1]{#1}%
%


\appendix
\section{Inner Product: Boundary-Adapted Basis}\label{app:ip}
In order to leverage the discrete Legendre transforms over 
Legendre-Gauss-Lobatto nodes, we favoured storage of space dependent variables 
in terms of their Legendre coefficients. As a result of this choice, all operations 
involving the boundary-adapted basis must be implemented via the Legendre transform.

Let us recall the definitions: $L_n(y)$ denote the Legendre polynomial of degree 
$n$ and the index set $\{0,1,\ldots,N-2\}$ is denoted by $\field{J}_1$. 
The boundary-adapted basis polynomials are given by
\begin{equation}
\begin{split}
&\phi_{p}(y_1)=L_p(y_1)+b_pL_{p+2}(y_1),\\
&b_p=-\frac{\kappa+\frac{1}{2}p(p+1)}{\kappa+\frac{1}{2}(p+2)(p+3)},\;p\in\field{J}_1,
\end{split}
\end{equation}
where $\kappa\in\field{C}$ is specific to the method chosen for temporal discretization
of boundary maps. The Fourier basis functions are given by $\psi_k(y)=\exp(iky)$ for 
$k\in\field{J}_2$. We can expand any arbitrary function in 
$\fs{L}^2(\field{I}\times\pi\field{I})$ in terms of the product of the two basis 
functions as follows:
\begin{equation}
\begin{split}
&f(\vv{y})=\sum_{p'=0}^{N-2}\sum_{q'\in\field{J}_2}
\hat{f}_{p'q'}\phi_{p'}(y_1)\psi_{q'}(y_2)\\
&+\sum_{q'\in\field{J}_2}\left[\eta_{0,q'}L_0(x_1)\psi_{q'}(y_2)
+\eta_{1,q'}L_1(x_1)\psi_{q'}(y_2)\right]\\
&=\sum_{p'=0}^N\sum_{q'\in\field{J}_2}\tilde{f}_{p'q'}L_{p'}(y_1)\psi_{q'}(y_2).
\end{split}
\end{equation}
Equating the coefficients in Legendre-Fourier basis (denoted by $\tilde{f}$) and boundary-adapted
Fourier basis (denoted by $\hat{f}$), we get
\begin{equation}
\left\{\begin{aligned}
&\tilde{f}_{0,q}=\hat{f}_{0,q}+\eta_{0,q}, \tilde{f}_{1,q}=\hat{f}_{1,q}+\eta_{1,q},\\
&\tilde{f}_{p,q}=\hat{f}_{p,q}+b_{p-2}\hat{f}_{p-2,q},\quad p=2,3,\ldots N-2,\\
&\tilde{f}_{N-1,q}=b_{N-3}\hat{f}_{N-3,q},\tilde{f}_{N,q}=b_{N-2}\hat{f}_{N-2,q}.
\end{aligned}\right.
\end{equation}
The linear system can be written as
\begin{equation}
\wtilde{F}=B\what{F} +\vv{e}_0\otimes\vs{\eta}_0 +\vv{e}_1\otimes\vs{\eta}_1
\end{equation}
where $B\in\field{C}^{(N+1)\times (N-1)}$ is given by
\begin{equation}
B=
\begin{pmatrix}
1   &    &      &      &      &\\
0   &  1 &      &      &      &\\
b_0 &  0 &  1   &      &      &\\
    & b_1&  0   &     1&      &\\
    &    &\ddots&\ddots&\ddots&\\
    &    &      &b_{N-4}&    0&1\\
    &    &      &       &b_{N-3}&0\\
    &    &      &       &     &b_{N-2}
\end{pmatrix},
\end{equation}
Let the normalization factors be denoted by 
$\gamma_k=\|L_k\|^2_2=2/(2k+1),\,k\in\field{N}_0,$ and introduce the diagonal matrix 
$\Gamma=\diag(\gamma_0,\gamma_1,\ldots,\gamma_{N})$, then inner products of the form 
$g_{pq}=\left(f,\phi_p\psi^*_q\right)$ can be understood as
\begin{equation}
\begin{split}
g_{pq}
&=\left(f,\phi_p\psi^*_q\right)\\
&=2\pi\left[\gamma_p\tilde{f}_{p,q}+b_p\gamma_{p+2}\tilde{f}_{p+2,q}\right],
\quad p\in\field{J},\,q\in\field{J}_2.
\end{split}
\end{equation}
The relationship above can be compactly stated as 
\begin{equation}
G=(g_{p,q})
=(g_{0,q},g_{0,q},\ldots,g_{N-2,q})^{\tp}
=Q\Gamma\wtilde{\vv{f}}M_2,
\end{equation} 
where $Q=B^{\tp}\in\field{C}^{(N-1)\times (N+1)}$ is referred to as the 
\emph{quadrature matrix}.

For CP and HF methods, the boundary-adapted polynomials are coupled with Fourier
components as:
\begin{equation}
\begin{split}
&\phi_{p,q}(y_1)=L_p(y_1)+b_{pq}L_{p+2}(y_1),\\
&b_{pq}=-\frac{\alpha_1\varpi_q+\frac{1}{2}p(p+1)}{\alpha_1\varpi_q+\frac{1}{2}(p+2)(p+3)},
\quad p\in\field{J}_1,\,q\in\field{J}_2.
\end{split}
\end{equation}
We can expand any arbitrary function in $\fs{L}^2(\field{I})$ in the 
basis as follows:
\begin{multline}
\wtilde{f}_q(y_1)
=\sum_{p=0}^{N-2}
\what{\wtilde{f}}_{p,q}\phi_{pq}(y_1)
+\eta_{0,q}L_0(x_1)\\+\eta_{1,q}L_1(x_1)
=\sum_{p=0}^N\wtilde{\wtilde{f}}_{p,q}L_{p}(y_1).
\end{multline}
Equating the coefficients in Legendre basis (denoted by $\wtilde{\wtilde{f}_q}$) and 
boundary-adapted basis (denoted by $\what{\wtilde{f}}$), the linear system can be written as
\begin{equation}
\wtilde{\wtilde{\vv{F}}}_q=B_{q,1}\what{\wtilde{\vv{F}}}_q
+\vv{e}_0\eta_{0,q}
+\vv{e}_1\eta_{1,q}
\end{equation}
where $\wtilde{B}_{q,1}\in\field{C}^{(N+1)\times (N-1)}$ is given by
\begin{equation}
\wtilde{B}_{q,1}=
\begin{pmatrix}
1   &    &      &      &      &\\
0   &  1 &      &      &      &\\
b_{0,q} &  0 &  1   &      &      &\\
    & b_{1,q}&  0   &     1&      &\\
    &    &\ddots&\ddots&\ddots&\\
    &    &      &b_{N-4,q}&    0&1\\
    &    &      &       &b_{N-3,q}&0\\
    &    &      &       &     &b_{N-2,q}
\end{pmatrix}.
\end{equation}
The inner products of the form 
$g_{pq}=\left(f,\phi_{pq}\right)$ can be understood as
\begin{equation}
\begin{split}
g_{pq}
&=\left(f,\phi_{pq}\right)\\
&=\left[\gamma_p\tilde{f}_{p,q}+b_{p,q}\gamma_{p+2}\tilde{f}_{p+2,q}\right],
\quad p\in\field{J},\,q\in\field{J}_2.
\end{split}
\end{equation}
The relationship above can be compactly stated as 
\begin{equation}
G=(g_{p,q})
=(g_{0,q},g_{0,q},\ldots,g_{N-2,q})^{\tp}
=\wtilde{Q}_{q,1}\Gamma\wtilde{\vv{f}},
\end{equation} 
where $\wtilde{Q}_{q,1}=\wtilde{B}_{q,1}^{\tp}\in\field{C}^{(N-1)\times (N+1)}$ is referred to as the 
quadrature matrix.
\end{document}